            \def\version{November 29th, 2008}     %
\documentclass[11pt,leqno]{amsart}
\usepackage{graphicx}
\usepackage{amsmath, amsthm, a4, latexsym, amssymb}
\usepackage{enumerate}
\def\@rmrk#1#2{\refstepcounter
    {#1}\@ifnextchar[{\@yrmrk{#1}{#2}}{\@xrmrk{#1}{#2}}}
 
\makeatletter\@addtoreset{equation}{section}\makeatother

\newfont{\bfit}{cmbxti10 scaled 2000}
\newfont{\biggi}{cmr12 scaled 2000}
 
 \newcommand{\eps}{\varepsilon}
 
 \newcommand{\Var}{{\rm Var}}

 \newcommand{\R}{\mathbb{R}}
 \newcommand{\Z}{\mathbb{Z}}
 \newcommand{\N}{\mathbb{N}}
 \newcommand{\ball}{{\mathcal B}}

 \newcommand{\prob}{\mathbb{P}}

 \newcommand{\me}{\mathbb{E}}

 \newcommand{\E}{\mathbb{E}}
 \newcommand{\1}{{\sf 1}}

 \newcommand{\FW}{{\mathfrak{W}}}

\newcommand{\To}[1]{\,\stackrel{#1}{\longrightarrow}\,}
 \newcommand{\Toi}[1]{\To{#1 \rightarrow \infty}}

 \def\CN{\mathcal{N}}

\newenvironment{Proof}[1]
{\vskip0.1cm\noindent{\sc #1}}{\vspace{0.15cm}}
\renewcommand{\subsection}{\secdef \subsct\sbsect}
\newcommand{\subsct}[2][default]{\refstepcounter{subsection}
\vspace{0.15cm}
{\flushleft\bf \arabic{section}.\arabic{subsection}~\bf #1  }
\nopagebreak\nopagebreak}
\newcommand{\sbsect}[1]{\vspace{0.1cm}\noindent
{\bf #1}\vspace{0.1cm}}
{\nopagebreak {\hfill{$\diamond$}}\\ }
\newtheorem{theorem}{Theorem}
\newtheorem{lemma}[theorem]{Lemma}

\newtheoremstyle{thm}{1.5ex}{1.5ex}{\itshape\rmfamily}{}
{\bfseries\rmfamily}{}{2ex}{}
\newtheoremstyle{rem}{1.3ex}{1.3ex}{\rmfamily}{}
{\itshape\rmfamily}{}{1.5ex}{}
\theoremstyle{rem}

\refstepcounter{subsubsection}
\def\thebibliography#1{\section*{Bibliography}
  \list%
  {\arabic{enumi}.}
    {\settowidth\labelwidth{[#1]}\leftmargin\labelwidth
    \advance\leftmargin\labelsep
    \parsep0pt\itemsep0pt
    \usecounter{enumi}}
    \def\newblock{\hskip .11em plus .33em minus .07em}
    \sloppy                   
    \sfcode`\.=1000\relax}

\def\CL{\mathcal{L}}

\def\FB{\mathfrak{B}}

\newcommand{\equ}[1]{(\ref{#1})}

\begin{document}
\title[Brownian intersection exponents]{Multiple intersection exponents}
\author[Klenke]{Achim Klenke}
\author[M{\"o}rters]{Peter M{\"o}rters}
\renewcommand{\thefootnote}{}
\subjclass[2000]{60J65.}
\keywords{Brownian motion, intersection exponent.}
\date{\version}
\begin{abstract}{\small Let $p\ge2$, $n_1\le\cdots\le n_p$ be positive integers and
$B_1^1, \ldots, B_{n_1}^1; \ldots; B_1^p, \ldots, B_{n_p}^{p}$ be
independent planar Brownian motions started uniformly on the boundary of
the unit circle. We define a $p$-fold intersection exponent $\varsigma(n_1,\ldots, n_p)$,
as the exponential rate of decay of the probability that the packets
$\bigcup_{j=1}^{n_i} B_j^i[0,t^2]$, $i=1,\ldots,p$, have no joint intersection.
The case $p=2$ is well-known and, following two decades of numerical and
mathematical activity,  Lawler, Schramm and Werner (2001) rigorously identified
precise values for these exponents. The exponents have not been investigated so far
for $p>2$. We present an extensive mathematical and numerical study, leading to
an exact formula in the case $n_1=1$, $n_2=2$, and several interesting
conjectures for other cases.}
\end{abstract}
\maketitle
\section{Introduction}

\subsection{Motivation and overview}

Finding \emph{exponents}, which describe the decay of some probabilities, and \emph{dimensions} of
some sets associated with stochastic models of physical systems is one of the core activities in
statistical physics. While in general one often has to resort to numerical methods to get a handle
on the values of the exponents, for planar models conformal invariance may help to answer these questions
explicitly, and there is now a substantial body of rigorous and non-rigorous methods available.
For example, by making the assumption that critical planar percolation behaves in a
conformally invariant way in the scaling limit and using ideas involving conformal field theory,
Cardy~\cite{Ca92} determined the asymptotic probability, as $N\to\infty$, that there exists
a two-dimensional critical percolation cluster crossing a rectangle.
A rigorous proof of Cardy's formula was later given by Smirnov~\cite{Sm01}. Following
considerable numerical work, see for example~\cite{LR78, Vo84} and references therein,
Saleur and Duplantier~\cite{SD87} predicted the
fractal dimension of the hull of a large percolation cluster using a non-rigorous Coulomb
gas technique. Rigorous versions of this result have been given based on Cardy's formula,
for example by Camia and Newman~\cite{CN06, CN07}.%
\smallskip

In~\cite{DK88} Duplantier and Kwon suggested that ideas of conformal field theory can also be used to
predict the probability of pairwise non-intersection between planar Brownian paths. Early research
by Burdzy, Lawler and Polaski~\cite{BLP89} and Li and Sokal~\cite{LS90} was of numerical nature, but
ten years later, Duplantier~\cite{Du98} gave a derivation based on non-rigorous
methods of quantum gravity, and soon after that Lawler, Schramm and Werner~\cite{LS01a, LS01b, LS02}
gave a rigorous proof based on the Schramm-Loewner evolution~(SLE), one of the greatest achievements in
probability in recent years. We also mention here some very recent developments with the long term aim of
making the quantum gravity approach rigorous, see Duplantier and Sheffield~\cite{DS08}, and
Rhodes and Vargas~\cite{RV08}.
\bigskip
\pagebreak[3]

In this paper we look at joint intersections of three or more planar Brownian paths, a question which
has been neglected so far in the literature, but which came up in our recent investigation of the
multifractality of intersection local times~\cite{KM05}.  In the simplest case, given three independent
Brownian paths $B^1$, $B^2$, $B^3$  started uniformly on the unit circle, we are interested in the asymptotic
behaviour, as $t\to\infty$, of the non-intersection probability
$$\prob\big\{ B^1[0,t] \cap  B^2[0,t] \cap B^3[0,t] = \emptyset \big\}.$$
Observe that this probability goes to zero, for $t\uparrow\infty$, as three, or any finite
number, of Brownian paths in the plane eventually intersect, see e.g.~\cite[Chapter 9.1]{MP09}.
Recall for comparison, that the non-intersection
exponents for three Brownian paths studied in the aforementioned papers deal with pairwise non-intersections,
i.e.~in the case of three Brownian motions either with  $$\prob\big\{ B^1[0,t] \cap  B^2[0,t]  = \emptyset,\,
B^2[0,t] \cap  B^3[0,t]  = \emptyset,\,
B^1[0,t] \cap  B^3[0,t]  = \emptyset \big\}, \quad \mbox{ or with }$$
$$\prob\big\{ B^1[0,t] \cap  \big( B^2[0,t] \cup B^3[0,t] \big)   = \emptyset \big\}.$$
Our study starts with the observation that, for positive integers $n_1,\ldots, n_p$
and independent planar Brownian motions
$$B_1^1, \ldots, B_{n_1}^1; \ldots; B_1^p, \ldots, B_{n_p}^{p},$$
nontrivial exponents
$$\varsigma(n_1, \ldots, n_p) = - \lim_{t\to\infty} \frac{2}{\log t}\,
\log \prob\Big\{ \bigcup_{j=1}^{n_1} B_j^1[0,t] \cap \ldots \cap \bigcup_{j=1}^{n_p} B_j^p[0,t]
= \emptyset \Big\}$$
exist, see Theorem~\ref{interdefi} and the subsequent remark.
In Theorem~\ref{res} we show that, for $2\le n_3\le\cdots\le n_p$, we have
$$\varsigma(1,2,n_3,\ldots, n_p)=2.\\[1mm]$$
These are the only exponents we could determine exactly beyond the well-known
case of $p=2$. Rigorous proofs of both theorems are given in Section~\ref{se.res}.
\smallskip

The bulk of this paper is devoted to the presentation of a detailed numerical study of the
values of the, in our opinion,  most interesting remaining exponents, see Section~\ref{se.num}.
One of the motivations of this study was to test the conjecture, motivated by Theorem~\ref{res}, that the value of the
exponents~$\varsigma(n_1,n_2,n_3,\ldots, n_p)$  depend only on the two smallest parameters.
This conjecture was not supported by our numerical investigations.
\smallskip

Finally, we remark that we have not been able to use either SLE techniques or quantum gravity
to derive even a non-rigorous exact prediction of the exponents if $p>2$. We hope however that our
numerical study triggers interest in this problem and that, as in the motivational examples discussed above,
future research will address the question of exact formulas for multiple intersection exponents.%
\bigskip

\subsection{Statement of the main theorems}

Let $p\ge 2$ and $n_1,\ldots, n_p$ be positive integers and $B_1^1, \ldots,
B_{n_1}^1; \ldots; B_1^p, \ldots, B_{n_p}^{p}$ independent planar
Brownian motions started uniformly on the unit circle~$\partial\ball(0,1)$. We define $p$ packets by
$$\FB^1(r):= \bigcup_{j=1}^{n_1} B_j^1\big[0,\tau^1_j(r)\big], \,\ldots\,,
\FB^p(r):=\bigcup_{j=1}^{n_p} B^p_j\big[0,\tau^p_j(r)\big],$$
where $\tau^i_j(r):=\inf\{t\ge 0 \colon |B_j^i(t)|=r\}$ and $r\ge 1$.

\begin{theorem}\label{interdefi}
The limit
$$\varsigma(n_1, \ldots, n_p):=
-\lim_{r\to\infty} \frac1{-\log r} \, \log \prob\Big\{ \FB^1(r) \cap \ldots \cap \FB^p(r)
= \emptyset \Big\}$$
exists and is positive and finite.
\end{theorem}

{\bf Remarks:}\ \\[-4mm]
\begin{itemize}
\item Using a standard argument, see \cite[Lemma 3.14]{La96}, one can replace the paths stopped upon
hitting the circle of radius~$r$, by paths running for $t=r^2$ time units. This leads to
the characterisation of the exponents given in the overview.\\[-1mm]
\item For $p=2$ all exponents are known, see \cite{LS01a, LS01b, LS02}. The technique used
to identify the exponents, which is based on the Schramm-Loewner evolution (SLE), does
not seem to allow us to identify the exponents for $p>2$. \\[-1mm]
\item We conjecture that one can strengthen this result, as this was done for~$p=2$ in~\cite{La95},
and show that there exists a constant $c>0$, depending on the starting points, such that
$$\lim_{r\to\infty} r^{\varsigma(n_1, \ldots, n_p)} \, \prob\Big\{ \FB^1(r) \cap \ldots \cap \FB^p(r)
= \emptyset \Big\} =c.$$
However, this is quite subtle and would go beyond the scope of this paper.
\end{itemize}
\smallskip

There is a trivial symmetry of the exponents, namely for every permutation $\sigma\in{\rm Sym}(p)$,
we have
$$\varsigma\big(n_1,\ldots, n_p\big)=\varsigma\big(n_{\sigma(1)},\ldots, n_{\sigma(p)}\big).$$
Moreover, there are two trivial monotonicity rules for these exponents
\begin{itemize}
\item[(A)] $\varsigma(n_1,\ldots, n_p) \le \varsigma(n_1,\ldots, n_{p-1})$, \\[-1mm]
\item[(B)] $\varsigma(n_1,\ldots, n_p) \le \varsigma(m_1,\ldots, m_{p})$, if $n_i \le m_i$ for $i=1,\ldots, p$.
\end{itemize}

As a result of the symmetry of the exponents, we may henceforth assume that the arguments of the exponents
are increasing in size, i.e. $n_1 \le \cdots \le n_p$. There is one interesting situation in which we can
determine the exponents explicitly.
\smallskip

\begin{theorem}\label{res}
We have $\varsigma(1,2,n_3,\ldots, n_p)=2$ for any $p\ge 2$ and $2\le n_3\le\cdots\le n_p$.
\end{theorem}

Note that to show this, by the monotonicity rules, it suffices to show that
\begin{equation}\varsigma(1,2,\,\,\,\stackrel{p-1}{\cdots}\,\,\,,2)=2\, .
\label{result}\end{equation}
The proof of this fact is based on the technique of hitting the intersection of $p-1$ Brownian paths by a further path,
using an idea of Lawler, see~\cite{La89} or~\cite[Section~3]{La91}, originally used to determine the exponent~$\varsigma(1,2)=2$.
\bigskip

\subsection{Conjectures}

In this section we formulate the main conjecture motivated by our numerical studies.
A detailed description of these studies and their outcomes will be given in Section~\ref{se.num}.

Let $p\in\N$ and $n_1,\ldots,n_p\in\N$ with $n_1\leq n_2\leq\ldots\leq n_p$. Define
$$k:=\min\big\{\ell\in\{2,\ldots,p\} \colon \,n_{\ell+1}>n_{\ell}\big\},$$
with $k:=p$ if the set is empty. We conjecture that
\begin{equation}
\label{EConj}
\varsigma(n_1,\ldots,n_p)=\varsigma(n_1,\ldots,n_k).
\end{equation}

In fact, this holds, by Theorem~\ref{res} for the case $k=2$, $n_k=2$,
and we have numerical evidence for
\begin{itemize}
\item
$\varsigma(1,1,2)=1.2503\pm0.0011$ to be compared with $\varsigma(1,1)=\frac54$\\[-1mm]
\item
$\varsigma(1,1,1,2)=1.02\pm0.004$ to be compared with
$\varsigma(1,1,1)=1.027\pm0.005$\\[-1mm]
\item
$\varsigma(1,3,3)=2.688\pm0.01$ to be compared with
$\varsigma(1,3)=\frac{13+\sqrt{73}}{8}=2.693000\ldots$\\[-1mm]
\item
$\varsigma(2,2,3)=2.937\pm0.01$ to be compared with
$\varsigma(2,2)=\frac{35}{12}=2.91666\ldots$\\[-1mm]
\item
$\varsigma(2,3,3)=3.767\pm0.06$ to be compared with
$\varsigma(2,3)=\frac{47+5\sqrt{73}}{24}=3.738334113\ldots$
\end{itemize}
\bigskip

\section{Proofs of Theorems~\ref{interdefi} and~\ref{res}.}\label{se.res}

\subsection{Proof of Theorem~\ref{interdefi}}

Denote by $x=(x_1^1, \ldots, x_{n_1}^1; \ldots; x_1^p, \ldots, x_{n_p}^{p})$ vectors
with $n_1+\cdots+n_p$ entries in $\R^2$, playing the role of configurations of our motions
at time zero. Consider
$$a_r:=\sup_{|x^i_j|=1} \prob_{x} \big\{ \FB^1(r) \cap \cdots \cap \FB^p(r)=\emptyset\big\},$$
where the subindex of $\prob$ indicates the starting points of the Brownian motions.
Using the strong Markov property and Brownian scaling, we get, for any $r,s\ge 1$,
\begin{align*}
a_{rs} & \le \sup_{|x^i_j|=1} \prob_{x}
\Big\{ \bigcup_{j=1}^{n_1} B^1_j\big[0,\tau^1_j(r)] \cap \cdots \cap \bigcup_{j=1}^{n_p} B^p_j\big[0,\tau^p_j(r)]
=\emptyset,\\ & \qquad\qquad\quad\;\;\bigcup_{j=1}^{n_1} B^1_j\big[\tau^1_j(r),\tau^1_j(rs)]
\cap \cdots\cap \bigcup_{j=1}^{n_p} B^p_j\big[\tau^p_j(r),\tau^p_j(rs)] =\emptyset\Big\}\\
& = \sup_{|x^i_j|=1} \me_{x}
\Big[ \1\Big\{\bigcup_{j=1}^{n_1} B^1_j\big[0,\tau^1_j(r)] \cap \cdots \cap \bigcup_{j=1}^{n_p} B^p_j\big[0,\tau^p_j(r)]
=\emptyset\Big\} \\ &
\qquad\qquad\qquad \times\prob_{(B^i_j(\tau^i_j(r))} \Big\{\bigcup_{j=1}^{n_1} B^1_j\big[\tau^1_j(r),\tau^1_j(rs)]
\cap \cdots\cap \bigcup_{j=1}^{n_p} B^p_j\big[\tau^p_j(r),\tau^p_j(rs)] =\emptyset\Big\} \Big]\\
& \le a_r a_s\, .
\end{align*}
Hence the function given by $b_t:=\log a_{2^t}$ is subadditive and,
by the subadditivity lemma, see e.g.~\cite[Lemma~5.2.1]{La91}, we thus have
$\lim_{t\to\infty} b_t/t = \inf_{t>0} b_t/t.$ Therefore,
$$\tilde\varsigma(n_1,\ldots,n_p):=
- \lim_{r\to\infty} \frac1{\log r} \log \sup_{|x^i_j|=1} \prob_{x} \big\{ \FB^1(r) \cap \cdots \cap \FB^p(r)=\emptyset\big\}$$
exists, and is positive.
\smallskip

Next, we show that we can replace the optimised starting points by starting points uniformly chosen from
the unit circle. Clearly, we have
\begin{equation}\label{no1}
\prob \big\{ \FB^1(r) \cap \cdots \cap \FB^p(r)=\emptyset\big\} \le
\sup_{|x^i_j|=1} \prob_{x} \big\{ \FB^1(r) \cap \cdots \cap \FB^p(r)=\emptyset\big\},
\end{equation}
where $\prob$ refers to the original scenario of Brownian motions started uniformly on the unit circle.
\smallskip

Conversely, using the Markov property, for $r>2$, we have
$$\begin{aligned}
\sup_{|x^i_j|=1} & \prob_{x} \big\{ \FB^1(r) \cap \cdots \cap \FB^p(r)=\emptyset\big\} \\
& \le
\sup_{|x^i_j|=1} \me_{x} \Big[ \prob_{(B^i_j(\tau^i_j(2)))}\Big\{\bigcup_{j=1}^{n_1} B^1_j\big[\tau^1_j(2),\tau^1_j(r)]  \cap \cdots \cap
\bigcup_{j=1}^{n_p} B^p_j\big[\tau^p_j(2),\tau^p_j(r)]=\emptyset\Big\} \Big].
\end{aligned}$$
By the Harnack principle, the law of the vector $(B^i_j(\tau^i_j(2)))$ is bounded, uniformly in~$x$, by a constant multiple of the uniform distribution on the $(n_1+\cdots+n_p)$-fold cartesian power of the circle $\partial\ball(0,2)$.
Denoting this constant by~$C$ and using Brownian scaling,
\begin{equation}\label{no2}
\begin{aligned}
\prob\Big\{\bigcup_{j=1}^{n_1} B^1_j\big[0,\tau^1_j(r/2)]  & \,  \cap \cdots \cap
\bigcup_{j=1}^{n_p} B^p_j\big[0,\tau^p_j(r/2)]=\emptyset\Big\}\\
&  \ge C^{-1} \,
\sup_{|x^i_j|=1} \prob_{x} \big\{ \FB^1(r) \cap \cdots \cap \FB^p(r)=\emptyset\big\} .
\end{aligned}
\end{equation}
Combining \eqref{no1} and \eqref{no2} yields that
$$\varsigma(n_1,\ldots,n_p):=
- \lim_{r\to\infty} \frac1{\log r} \log \prob\big\{ \FB^1(r) \cap \cdots \cap \FB^p(r)=\emptyset\big\}$$
exists and coincides with $\tilde\varsigma(n_1,\ldots,n_p)$. Note, finally, that the monotonicity rule (A) implies
that $\varsigma(n_1,\ldots,n_p)\le \varsigma(n_1,n_2)<\infty$, and hence the exponents are positive and finite.
{\nopagebreak {\hfill{$\diamond$}}\\ }
\smallskip

\subsection{Proof of Theorem~\ref{res}}

Recall that it suffices to show~\eqref{result}. We start by formulating the key lemma.
We let $W^1, \ldots, W^p$ be independent Brownian paths. For $r,s>0$ denote by $\tau^i(x,r)$
the first hitting time by the motion $W^i$ of the circle $\partial
\ball(x,r)$ with centre $x$ and radius $r$, and let $\tau^i(x,r,s)$
be the first hitting time of $\partial \ball(x,s)$
after~$\tau^i(x,r)$.
\begin{lemma}\label{revISE}
Fix $x\in\ball(0,1)$. Suppose that $W^1, \ldots, W^p$ are
independent Brownian paths started uniformly on the circle
$\partial\ball(0,2)$. Define the set
\begin{equation}
\label{EL2.1}
\FW:=\bigcap_{j=2}^p W^j[0, \tau^j(0,4)]
\end{equation}

and the events

\begin{equation}
\label{EL2.2}
\begin{aligned}
E_{x,r}&=\big\{W^1[0,\tau^1(x,r/2)]
\cap \FW=\emptyset\big\},\\
N_{x,r}&=\big\{W^1[0,\tau^1(x,r/2,r)]
\cap\FW\neq\emptyset\big\},\\
H_{x,r}&=\big\{\tau^i(x,r/2)<\tau^i(0,4)\mbox{ \,for all\, }
i=1,\ldots,p\big\}.
\end{aligned}
\end{equation}
Then
$$\begin{aligned}
&\liminf_{r\downarrow 0}\frac 1{|\log r|} \,
\log \prob\big[E_{x,r}\cap N_{x,r}\,\big|\,H_{x,r}\big]
\ge - \varsigma(1,2,\stackrel{p-1}{\cdots},2).
\end{aligned}$$
\end{lemma}
\medskip

Let us first see how~\eqref{result} follows from this lemma.
Let
$$\tau=\inf\big\{ t>0 \, : \, W^1(t)\in \FW \big\}.$$
Now let $\mathfrak B$ be a collection of pairwise disjoint discs
of fixed radius $0<r<1/2$ with centres in the disc
$\ball(0,1)$, which has cardinality at least $(2r)^{-2}$. Then,
obviously,
$$\begin{aligned}
1  \;\ge\; \prob\big\{ W^1[0,\tau^1(0,4)]
\cap\FW\neq\emptyset \big\}
\;\ge\; \sum_{\ball \in \mathfrak B} \prob\big\{ W^1(\tau) \in \ball, \, \tau<\tau^1(0,4) \big\}\, .
\end{aligned}$$
Now, fix a disc $\ball=\ball(x,r)\in \mathfrak B$. The event $\{ W^1(\tau) \in \ball, \, \tau<\tau^1(0,4) \}$ is implied by the events
$$
\begin{aligned}
E_{x,r}\cap N_{x,r} \cap \{\tau^1(x,r/2)<\tau^1(0,4)\}.
\end{aligned}
$$
Recall that $$\prob\big[H_{x,r}\big]=\prob\big\{
\tau^1(x,r/2)<\tau^1(0,4)\big\}^p = r^{o(1)}\, .$$ Combining this
with Lemma~\ref{revISE}, for any $\eps>0$ and sufficiently small
$r>0$,
$$\prob\big\{ W^1(\tau) \in \ball, \, \tau<\tau^1(0,4) \big\}
\ge r^{\varsigma(1,2,\ldots,2)+\eps}.$$ This implies
$$\begin{aligned} 1 & \ge \sum_{\ball \in \mathfrak B} r^{\varsigma(1,2,\ldots,2)+\eps}
\ge r^{-2+\varsigma(1,2,\ldots,2)+2\eps}\, ,\end{aligned}$$
and therefore $\varsigma(1,2,\ldots,2)\ge 2-2\eps$. The lower bound follows
as $\eps>0$ was arbitrary, and the upper bound in \eqref{result} follows from
$\varsigma(1,2,\ldots,2)\le \varsigma(1,2)=2$, as
is known from \cite{La89, La91}.{\nopagebreak {\hfill{$\diamond$}}\\ }
\medskip
\pagebreak[2]

\begin{Proof}{Proof of Lemma~\ref{revISE}.}
Before we describe the technical details we sketch the idea of
the proof. Since the paths of $p$ planar Brownian motions
intersect with positive probability, by Brownian scaling,
the conditional probability of  $N_{x,r}$ given $H_{x,r}$ is
bounded from below as $r\to0$. Hence this
condition can be neglected when computing the probability in
Lemma~\ref{revISE}. For $j=1,\ldots,p$ we decompose the paths
$W^j$ into the pieces $W^j[0,\tau^j(x,r/2)]$ and $W^j[\tau^j(x,r/2),\tau^j(0,4)]$. By
time reversal for $W^j[0,\tau^j(x,r/2)]$, we can
compare the probability in question with the non-intersection
probability for packets of size $n_1=1$, $n_2=\cdots=n_p=2$,
which is of order~$\approx r^{\varsigma(1,2,\ldots,2)}$.
\smallskip

We now come to the technical details, see the appendix in~\cite{MS09} for the necessary
facts about Brownian excursions between concentric spheres.  Let $\varrho^1=r$ and
$\varrho^j=r/2$ for $j=2,\ldots,p$. Conditioned on
$\{\tau^i(x, \varrho^j/2)<\tau^i(x,3)\}$ the path
$W^i[0,\tau^i(x,\varrho^j/2)]$ is contained in an excursion from
$\partial\ball(x,3)$ to $\partial\ball(x,\varrho^j/2)$. The time-reversal of this
excursion is contained in the path of a Brownian motion ${\widetilde W}^i$
started uniformly on $\partial\ball(x,\varrho^j/2)$ and stopped upon
reaching $\partial\ball(x,3)$, say at time $\widetilde\tau^i(x,3)$.
Analogously to \equ{EL2.1} and \equ{EL2.2} define the set
$$\widetilde\FW=\bigcap_{j=2}^p \big(
{\widetilde W}^j[0,\widetilde\tau^j(x,3)] \cup W^j[\tau^j(x,r/4,r/2),\tau^j(0,4)]\big),$$
and the events
$$\begin{aligned}
\widetilde E_{x,r}&=\big\{{\widetilde W}^1[0,\widetilde\tau^1(x,3)] \cap
\widetilde\FW=\emptyset\big\},\\
\widetilde N_{x,r}&=\bigg\{
\bigcap_{j=1}^p W^j[\tau^j(x,\rho^j/2), \tau^j(x,\rho^j/2,\rho^j)]
\not= \emptyset\bigg\},\\
\widetilde H_{x,r}&=\big\{ \tau^j(x,\varrho^j/2)<\tau^j(x,3)
\mbox{ \,for all\, }j=1,\ldots,p\big\}.\\[1mm]
\end{aligned}
$$
Note that
$W^1[0,\tau^1(x,\rho^1)]\cap \ball(x,r/2)=\emptyset$ and $W^j[\tau^j(x,\rho^j/2),\tau^j(x,\rho^j/2,\rho^j))\subset \ball(x,r/2)$ for $j=2,\ldots,p$. Hence
$$W^1[0,\tau^1(x,\rho^1)]\cap \big(\FW\setminus \widetilde\FW\big)
\subset W^1[0,\tau^1(x,\rho^1)]\cap\bigcap_{j=2}^p W^j[\tau^j(x,\rho^j/2),\tau^j(x,\rho^j/2,\rho^j))=\emptyset$$
which implies $\widetilde E_{x,r}\subset E_{x,r}$.
Note that trivially, we have $\widetilde H_{x,r}\subset H_{x,r}$ and $\widetilde N_{x,r}\subset N_{x,r}$ which implies
\begin{equation}
\label{EL3.3.1}
E_{x,r}\cap N_{x,r}\cap H_{x,r}
\;\supset\;
\widetilde E_{x,r}\cap \widetilde N_{x,r}\cap \widetilde H_{x,r}.
\end{equation}

Finally, note that
\begin{equation}
\label{EL3.3.2}
f(x,r):=\frac{\prob[\widetilde H_{x,r}]}{\prob[H_{x,r}]}=\frac{\prob\big\{\tau^1(x,\varrho^1/2)<\tau^1(x,3)\big\}^p}
{\prob\big\{\tau^1(x,r/2)<\tau^1(0,4)\big\}^p}\geq\frac12
\end{equation}
for all $x$ and for sufficiently small values of $r>0$.

By \equ{EL3.3.1}, \equ{EL3.3.2} and the definition of the conditional probability, we conclude
\begin{equation}
\label{EL2.4}
\prob\big[E_{x,r}\cap N_{x,r}\,\big|\,H_{x,r}\big]
\ge f(x,r)\,\prob\big[\widetilde E_{x,r}\cap \widetilde N_{x,r}\big|\widetilde
H_{x,r}\big].\\[1mm]
\end{equation}
Fix $\eps>0$. Invoking the definition of the exponent, the Harnack principle and Brownian scaling,
for sufficiently small $r>0$,
$$\prob\big[ \widetilde E_{x,r} \, \big| \, \widetilde H_{x,r}\big]
\geq r^{\varsigma(1,2,\ldots,2)+\eps}.$$

Define the compact sets
$$\begin{aligned}
C&:=\big\{y=(y^1,\ldots,y^p): y^j\in\partial\ball(0,\varrho^j/2)\mbox{\; for }j=1,\ldots,p\big\}\mbox{\; and}\\
D&:=
\big\{z=(z^1,\ldots,z^p): z^j\in\partial\ball(0,\varrho^j)\mbox{\; for }j=1,\ldots,p\big\}.
\end{aligned}
$$
For $y\in C$ and $z\in D$ let $(\bar W^j, \;j=1,\ldots,p)$ be an independent family of Brownian motions where each motion $\bar W^j$ is started at $y^j$ and is conditioned to leave $\ball(0,\varrho^j)$ at $z^j$ (at time $\bar\tau^j$). Denote by $\prob_{y,z}$ the corresponding probability measure. It is easy to see that the map
$$\phi:C\times D\to[0,1],\quad (y,z)\mapsto\prob_{y,z}\big\{\bar W^1[0,\tau^1]\cap\ldots \cap \bar W^p[0,\tau^p]\neq\emptyset\big\}$$
is continuous and strictly positive, and independent of~$r$ by Brownian scaling. Hence
$$c:=\inf_{y\in C,\, z\in D}\phi(y,z)>0.$$
We infer that
$$\prob\big[ \widetilde N_{x,r} \,  \big| \, \widetilde E_{x,r} \cap \widetilde H_{x,r}\big] \ge c>0.\\[1mm]
$$
Hence, combing our results, for sufficiently small~$r>0$
$$\prob\big[\widetilde E_{x,r}\cap \widetilde N_{x,r}\big|\widetilde
H_{x,r}\big] =\prob\big[ \widetilde E_{x,r} \, \big| \, \widetilde H_{x,r}\big]
\prob\big[ \widetilde N_{x,r} \,  \big| \, \widetilde E_{x,r} \cap \widetilde H_{x,r}\big]
\ge c\, r^{\varsigma(1,2,\ldots,2)+\eps},\\[1mm]$$
and this completes the proof as $\eps>0$ was arbitrary.
{\nopagebreak {\hfill{$\diamond$}}\\ }
\end{Proof}

\section{Simulations}\label{se.num}
\label{S2}
To get hold of those exponents which we could not determine explicitly,
we have performed Monte Carlo simulations. This has successfully generated conjectures in the
$p=2$ case, see Duplantier and Kwon~\cite{DK88}, Li and
Sokal~\cite{LS90} and Burdzy, Lawler and Polaski~\cite{BLP89}.

\subsection{The general scheme.}
\label{S2.1}
Before we list and analyse the simulated data, we explain how we
got it. Fix positive integers $p$ and $n_1,\ldots,n_p$. The aim is
to get an estimate on $\varsigma(n_1,\ldots,n_p)$. Instead of
Brownian motions we simulate two-dimensional symmetric nearest
neighbour random walks. As it reduces computing effort, we work
with boxes rather than with discs. (For comparison we have performed
some of the simulations also with discs and there was no significant
difference in the results.) First we fix an increasing
sequence of box half-lengths $L_0,\ldots,L_K$ (in most cases
$L_{k+1}=\lfloor 1.1\cdot L_k\rfloor$ and the maximal value
$m=L_{L}$
restricted to 20000, 40000 or 80000) and the sample size $N$ of the
simulation.

\textbf{Step 1.} We start $n_1+\ldots+n_p$ independent random
walks at the origin $0\in\Z^2$ and stop each of them when it hits
$\partial(\{-L_0,\ldots,L_0\}^2)$. This defines the starting
positions of the random walks.

\textbf{Step 2.} Assume we are at level $k$ (after Step 1 we are
at level $k=1$). Independently run the random walks until they hit
$\partial(\{-L_{k},\ldots,L_{k}\}^2)$. Separately, keep track of
the set $A_{k,i}\subset \{-L_{k}+1,\ldots,L_{k}-1\}^2$ of points
that are visited by the $i$th package of $n_i$ random walks
\emph{before} hitting $\partial(\{-L_k,\ldots,L_k\}^2)$ (after
Step 1).

If $A_{k,1}\cap\ldots\cap A_{k,p}=\emptyset$, then we say that we
have \emph{survived} level $k$ and we enter level $k+1$ (that is, we perform
Step 2 again with $k$ replaced by $k+1$). Otherwise we stop this
sample and start a new simulation in Step 1.

By $N_k$ we denote the number of samples that have survived level
$k$. Clearly, $N_0=N$. We should have
$$
N_k/N\approx (L_k/L_0)^{-\varsigma(n_1,\ldots,n_p)}.
$$
Hence in a double logarithmic plot of $\log(N_k)$ against
$\log(L_k)$ the points should be on a line with slope
$-\varsigma(n_1,\ldots,n_p)$. Linear regression then gives an
estimate for the exponent $\varsigma(n_1,\ldots,n_p)$.

As it turns out that a line can be fitted well only for large
values of $L_k$, we have neglected the small values of $L_k$ in
order to get a reasonable estimate for
$\varsigma(n_1,\ldots,n_p)$. In Figure ~\ref{Fe11-linreg} below we plotted
the data points used for the linear regression with solid
circles, the other points with hollow circles.

As can be seen from Figure~\ref{Fe11-linreg}, for $\xi(1,1)$ this gives a pretty good estimate of the exact value $\frac54$, even with a moderate computing effort of about 2000 hours CPU time. However, for $\xi(1,1,1)$ the points tend to lie on a straight line only for large values of $L_k$ and thus require
\begin{enumerate}[(i)]
\item
a large maximal box size $m=L_K$ and thus a big computer memory of size $(2m+1)^2$ bytes in order to keep track of the visited points,
\item
a large sample size $N_0$ in order that $N_K\approx N_0\cdot(L_K/L_0)^{-\xi(1,1,1)}$ is big enough to obtain reliable data from the simulation.
\end{enumerate}
Since the CPU time we need for each sample grows with $m$, (i) and (ii) imply that we need huge amounts of CPU time. Furthermore, with huge sample sizes and box sizes, we run into the order of the cycle length of the common 48 bit linear congruence random number generators.

The computations were performed on different computers, mainly on two parallel Linux clusters at the University of Mainz on Opteron 2218 processors with 2.6GHz and on Opteron 244 processors with 1.8GHz.
The programme code is written in C. As random number generator we used {\sf drand64()}, a 64 bit linear congruence generator following the rule
$$          r_{n+1} = (a r_n + c) \mod 2^{64}$$
with
$$a = 6364136223846793005\quad\mbox{and}\quad c=1$$
(see \cite[pp106-108]{Kn05}).

\def\kmin{{k_{\mathrm{min}}}}
\def\Lmin{{L_{\mathrm{min}}}}
\def\Lmax{{L_{\mathrm{max}}}}
The linear regression method does not give a quantitative estimate on the statistical error. In order to get such an error estimate we did the following. Having in mind that the systematic error is large for small box sizes, we choose a minimal box number $\kmin\in\{1,\ldots,K-1\}$ and neglect the data from all smaller boxes. Furthermore, we pretend that the asymptotics for $p_L$ is exact for $k\geq\kmin$, that is,
\begin{equation}
\label{ESS1a}
p_{L_k}=C\,L_k^{-\varsigma}\quad\mbox{for all }k\geq \kmin
\end{equation}
for some $C>0$. In particular, the conditional probability to have no multiple intersections before leaving $B_{L_{k+1}}$ given there is no multiple intersection before leaving $B_{L_{k}}$ is
$$\bar p_k:=\frac{p_{L_{k+1}}}{p_{L_k}}=\left(\frac{L_k}{L_{k+1}}\right)^{-\varsigma}
=: q_k^{-\varsigma}.
$$
Here the likelihood function for the observation $$(N_\kmin,N_{\kmin+1},\ldots,N_K)=n:=(n_\kmin,n_{\kmin+1},\ldots,n_K)$$ is
\begin{equation}
\label{ESS1b}
\begin{aligned}
L_n(\varsigma)
&=C(n)\prod_{l=\kmin}^{K-1}\bar p_l{}^{n_{l+1}}(1-\bar p_l)^{n_l-n_{l+1}}\\
&=C(n)\prod_{l=\kmin}^{K-1}\bar q_l{}^{\varsigma n_{l+1}}(1-\bar q_l{}^\varsigma)^{n_l-n_{l+1}}
\end{aligned}
\end{equation}
for some $C(n)>0$. The log-likelihood function is
\begin{equation}
\label{ESS1c}
\CL_n(\varsigma)=\log C(n)+\sum_{l=\kmin}^{K-1}
\Big(n_{l+1}\varsigma\,\log(q_l)\;
+\;(n_l-n_{l+1})\,\log\big(1-q_l^\varsigma\big)\Big).
\end{equation}
The maximum likelihood estimator (MLE) $\hat \varsigma$ is defined by
\begin{equation}
\label{ESS1d}
\CL_n\big(\hat\varsigma\big)=\sup_{\varsigma>0}\CL_n(\varsigma).
\end{equation}

We compute the derivatives
\begin{equation}
\label{ESS1e}
\CL_n'(\varsigma)=\sum_{l=\kmin}^{K-1}n_{l+1}\,\log(q_l)
\;-\;\sum_{l=\kmin}^{K-1}(n_l-n_{l+1})
\frac{\log(q_l)\,q_l^\varsigma}{1-q_l^\varsigma}
\end{equation}
and
\begin{equation}
\label{ESS1f}
\CL_n''(\varsigma)=-\sum_{l=\kmin}^{K-1}(n_l-n_{l+1})
\frac{(\log(q_l))^2\,q_l^\varsigma}{(1-q_l^\varsigma)^2}.
\end{equation}
Clearly, $\CL_n''(\varsigma)<0$, hence $\varsigma\mapsto\CL_n(\varsigma)$ is strictly concave and thus $\hat\varsigma$ is the unique solution of
\begin{equation}
\label{ESS1f2}
\CL_n'(\hat\varsigma)=0.
\end{equation}
Hence, for given data, the MLE can easily be computed numerically (we used a Newton approximation scheme).

Denote by $\hat\varsigma_{n_0}$ the MLE for sample size $n_0$.
By standard theory for MLEs, $(\hat\varsigma)_{n_0\in\N}$ is consistent and asymptotically normally distributed. In fact, by Corollary~6.2.1 of \cite{Le83},
\begin{equation}
\label{ESS2a}
\hat\varsigma_{n_0}\Toi{n_0}\varsigma\quad\mbox{stochastically.}
\end{equation}

Furthermore, by \cite[Corollary 6.2.3]{Le83}, $(\hat\varsigma)_{n_0}$ is asymptotically efficient (that is, optimal) and by \cite[Theorem 6.2.3]{Le83} (with $\CN_{0,1}$ the standard normal distribution)
\begin{equation}
\label{ESS2b}
\sqrt{n_0 I(\varsigma)}
\big(\hat\varsigma_{n_0}-\varsigma\big)\Toi{n_0}\CN_{0,1}\quad\mbox{in distribution}.
\end{equation}
Here
\begin{equation}
\label{ESS2c}
I(\varsigma)=-\E[\CL_N''(\varsigma)|N_0=1]=p_{L_\kmin}
\sum_{l=\kmin}^{K-1}\left(\prod_{m=\kmin}^{l-1}\bar p_m\right)\big(1-\bar p_l\big)\frac{(\log(q_l))^2\,q_l^\varsigma}{(1-q_l^\varsigma)^2}
\end{equation}
is the Fisher information for one sample. As we do not know the true value of $\varsigma$ and since we do not know $p_{L_\kmin}$, we replace $I(\varsigma)$ by
$$
I_n(\varsigma)=-\frac{1}{n_0}\CL_n''(\varsigma).
$$
By the law of large numbers $I_N(\varsigma)\Toi{n_0}I(\varsigma)$ almost surely, uniformly in $\varsigma$ in compact sets. Hence by \equ{ESS2a}, we have $I_N(\hat\varsigma)\Toi{n_0}I(\varsigma)$ stochastically.
Hence we use
\begin{equation}
\label{ESS2d}
\widehat\sigma{}\,^2:=-1/\CL''_N(\hat\varsigma)
\end{equation}
as an estimator for the variance of $\hat\varsigma$ and obtain
\begin{equation}
\label{ESS2e}
\frac{\hat\varsigma-\varsigma}{\widehat\sigma}\Toi{n_0}\CN_{0,1}
\quad\mbox{in distribution}.
\end{equation}
Concluding, an asymptotic 95\% confidence interval for $\varsigma$ is given by
\begin{equation}
\label{ESS1h}
\big[\hat\varsigma-2\,\widehat\sigma,\,\hat\varsigma+2\,\widehat\sigma\big].
\end{equation}

We have performed the simulations for the exponents
$\varsigma(1,1)$ and $\varsigma(2,2)$ as benchmark problems, and then did the simulations on a larger scale for
$$\varsigma(1,1,1),\quad\varsigma(1,1,2),\quad\varsigma(1,1,1,1),\quad\varsigma(1,1,1,2).
$$

\subsection{Two-level scheme.}
\label{S2.2}
The simulations turn out to be very time-consuming, especially for the exponents with a larger numerical value. In order to get a more efficient scheme in this situation consider the following simplification of the simulation scheme presented above:

Assume there are only three box sizes, $L_0$ (about 30), $L_1$ (about (10\,000) and $L_2=2L_1$. Then \equ{ESS1f2} can be solved explicitly and the maximum likelihood estimator for $\varsigma$ is
$$\hat\varsigma = -\frac{\log(n_2/n_1)}{\log(2)}.$$
In order to reduce the variance of $\hat \varsigma$ we have to increase $N_1$, that is the sample size $n_0$. However, since it takes much CPU time to obtain a sample that contributes to $N_1$, we may wish to use this very sample as the starting point for a number $m$ of trials running from box size $L_1$ to $L_2$. Assume that $x$ among these  $m$ trials have survived until $L_2$ (that is, have reached the boundary of the $L_2$-box without producing a multiple intersection), then $p_S=\frac{x}{m}$ is an estimator for the conditional probability of producing no multiple intersection until leaving the $L_2$-box for the given realisation $S$ of the paths of all walks in the $L_1$-box. Now we can prescribe the number $n=n_1$ of ``master samples'' and for $i=1,\ldots,n$ let $x_i$ be the corresponding number of surviving trials and write $\widehat {p_i}:=x_i/m$. Hence for
$$p:=\frac{p_{L_2}}{p_{L_1}}=\E[p_S]$$
we get the unbiased estimator
$$\hat p=\frac1n\sum_{l=1}^n\widehat {p_i}.$$
The unbiased estimator for the variance of $\hat p$ is
$$\widehat {\sigma_p^2}=\frac{1}{n(n-1)}\sum_{l=1}^n(\widehat{p_i}-\hat p)^2.$$
From $\hat p$ and $\widehat {\sigma_p^2}$ we obtain the estimators for $\varsigma$ and the variance $\sigma^2$ of $\hat\varsigma$
\begin{equation}
\label{EESS2}
\hat\varsigma =-\frac{\log(p)}{\log(2)}
\qquad\mbox{and}\qquad
\widehat{\sigma^2}=\frac{\widehat {\sigma_p^2}}{(\log(2)\,\hat p)^2}.
\end{equation}

We have employed this scheme for the exponents with numerical values larger than 2, and we explain now why it is more efficient in these cases.

The expected time planar random walk needs to go from $\partial(\{-L,\ldots,L\}^2)$ to $\partial(\{-L-1,\ldots,L+1\}^2)$ is of order $L$. The probability that a given sample ever reaches $\partial(\{-L-1,\ldots,L+1\}^2)$ is of order $L^{-\varsigma}$. Hence (if we stop the simulation as soon as the first multiple intersection is detected) the expected CPU time for each sample until box size $L_1$ is of order
$$\sum_{L=L_0}^{L_1}L^{1-\varsigma}.$$
For $\varsigma>2$ this sum is of order $1$, for $\varsigma\leq 2$, it is of order $L_1^{2-\varsigma}$. Now the probability that a sample reaches box size $L_1$ without producing a multiple intersection is of order $L_1^{-\varsigma}$. Hence the expected CPU time needed for simulating a ``master sample'' is of order $L^{2\vee \varsigma}$. On the other hand, each of the trials started from the master sample needs an expected CPU time of order $L_1^2$. Hence for $\varsigma>2$ we can run $m=L_1^{\varsigma-2}$ trials without increasing the CPU significantly.

In order to make a good choice for $m$, compute the variance of $\hat p$
$$\Var[\hat p]=n^{-1}\Var[p_S]+\frac{1}{mn}\E[p_S(1-p_S)]\leq
n^{-1}\Var[p_S]+\frac{1}{mn}\E[p_S].$$
The quantities $\Var[p_S]$ and $\E[p_S]\approx 2^{-\varsigma}$ can be estimated from a test simulation as well as the expected CPU time $T_1$ to produce a master sample and the expected time $T_2$ used for each subsequent trial.
Now it is an optimisation problem for the total CPU time $n(T_1+mT_2)$ versus the variance $\Var[\hat p]$. For some of the simulations we have done test runs and solved the optimisation problem. Here $m=1000$ turned out to be a reasonable choice that we have then used in all simulations.

We have performed the simulations according to this scheme with $L_0=30$, $L_1=10\,000$, $L_2=20\,000$ and $m=1000$ for the exponents
$$\varsigma(1,3,3),\quad\varsigma(2,2,2),\quad\varsigma(2,2,3),
\quad\varsigma(2,3,3),\quad\varsigma(2,2,2,2).$$
\clearpage

\subsection{Numerical results.}
\label{S2.3}
We present our estimated values $\hat\varsigma$ together with a statistical error of $2\sigma$. For the systematic error it is hard to make a good judgement. From the graphical representation of the results (see below) it seems that for $\varsigma(1,1,2)$ the systematic error is of a smaller order than the statistical error. For $\varsigma(1,1,2)$ and $\varsigma(1,1,1)$ it is presumably of the same order. Finally, for $\varsigma(1,1,1,1)$ and, even worse for $\varsigma(1,1,1,1,1)$ we seem to systematically underestimate the values. It would require a lot larger $\Lmax$ to get more accurate results. For that reason we have not taken too much effort to reduce the statistical error. However, we give the results of the simulations just to provide an idea of the possible values.

\begin{table}[h]
\label{TS1}
$$\def\arraystretch{1.2}
\begin{array}{|r|l|l|l|r|l|r|r|}
\hline
\mbox{exponent}&\multicolumn{1}{c|}{\hat\varsigma}&2\widehat\sigma&\mbox{rigorous}&\Lmin&\Lmax&n_0/10^6&\multicolumn{1}{c|}{\mbox{CPU}}\\
&&&&&&&\multicolumn{1}{c|}{\mbox{ time/h}}\\\hline
\varsigma(1,1)&1.2502&0.001&5/4&1069&20\,000&500&2\,064\\
\varsigma(2,2)&2.9188&0.0033&\frac{35}{12}=2.9167&163&20\,000&40\,000&1\,879\\
\varsigma(1,1,1)&1.027&0.005&[1/2,\,5/4]&18\,575&80\,000&60&8\,262\\
\varsigma(1,1,2)&1.2503&0.0011 &[1, 5/4]&1069&80\,000&200&5\,858\\
\varsigma(1,1,1,1)&0.877&0.006 &[1/4, 5/4]&39\,813&80\,000&20&18\,262\\
\varsigma(1,1,1,2)&1.02&0.004 &[1/2, 5/4]&27\,194&40\,000&200&35\,212\\
\varsigma(1,1,1,1,1)&0.74&0.02 &[1/8, 5/4]&27\,194&40\,000&0.74&1\,147\\
\hline
\end{array}
$$
\caption{Numerical results obtained from the first simulation scheme.}
\end{table}

\begin{table}[h]
\label{TS2}
$$\def\arraystretch{1.2}
\begin{array}{|r|l|l|l|r|r|}
\hline
\mbox{exponent}&\multicolumn{1}{c|}{\hat\varsigma}&2\widehat\sigma&\mbox{rigorous}&n&\multicolumn{1}{c|}{\mbox{CPU}}\\&&&&&\multicolumn{1}{c|}{\mbox{time/h}}\\
\hline
\varsigma(1,3,3)&   2.688&    0.01  & [2,\,(13+\sqrt{73})/8]          &18\,100&61\,860\\
\varsigma(2,2,2)&   2.786&    0.01  & [2,\,35/12] &16\,000&47\,943\\
\varsigma(2,2,3)&   2.937&    0.01  & [2,\,35/12] &23\,000&116\,888\\
\varsigma(2,3,3)&   3.767&    0.057 & [2,\,35/12] & 1\,000&179\,543\\
\varsigma(2,2,2,2)& 2.664&    0.01  & [2,\,35/12] &16\,000&63\,496\\
\hline
\end{array}
$$
\caption{Numerical results obtained from the second simulation scheme.}
\end{table}
\clearpage

\subsection{Detailed Data.}
\subsubsection{Exponent $\varsigma(1,1)$.}
The exact value $\varsigma(1,1)=5/4$ is known. This simulation is used as a benchmark test for our simulation.
\def\ftl{\begin{equation*}\tiny\setlength{\doublerulesep}{0pt}%
\def\arraystretch{1.1}\quad\begin{array}{||r|r||}\hline\hline L_k&n_k\\%
\hline}
\def\etl{\hline\hline\end{array}\quad\end{equation*}}
\def\ntl{\hline\hline\end{array}\quad\begin{array}{||r|r||}%
\hline\hline L_k&n_k\\\hline\hline}

\ftl
   30&  500000000\\
   33&  455164209\\
   36&  414185142\\
   39&  379373384\\
   42&  349383901\\
   46&  315390855\\
   50&  286840826\\
   55&  257075021\\
   60&  232385705\\
   66&  207870728\\
   72&  187620511\\
   79&  168084821\\
   86&  151902122\\
   94&  136553134\\
  103&  122326905\\
\ntl
  113&  109366745\\
  124&  97714439\\
  136&  87320799\\
  149&  78109962\\
  163&  69978568\\
  179&  62384176\\
  196&  55800459\\
  215&  49786852\\
  236&  44382636\\
  259&  39563995\\
  284&  35298660\\
  312&  31418279\\
  343&  27932867\\
  377&  24837149\\
  414&  22109889\\
\ntl
  455&  19660552\\
  500&  17483797\\
  550&  15525080\\
  605&  13788917\\
  665&  12253892\\
  731&  10890052\\
  804&  9669275\\
  884&  8589857\\
  972&  7631215\\
 1069&  6776772\\
 1175&  6020939\\
 1292&  5347118\\
 1421&  4747333\\
 1563&  4214131\\
 1719&  3741150\\
\ntl
 1890&  3323382\\
 2079&  2950258\\
 2286&  2620862\\
 2514&  2327160\\
 2765&  2066024\\
 3041&  1834523\\
 3345&  1628901\\
 3679&  1446024\\
 4046&  1283655\\
 4450&  1140213\\
 4895&  1012659\\
 5384&  898680\\
 5922&  797641\\
 6514&  708293\\
 7165&  628813\\
\ntl
 7881&  557957\\
 8669&  495180\\
 9535&  439662\\
10488&  389839\\
11536&  345918\\
12689&  307046\\
13957&  272420\\
15352&  241798\\
16887&  214746\\
18575&  190486\\
20000&  173506\\
&\\&\\&\\&\\\etl
\medskip

Values used for the fit: $L_k=1069\ldots20\,000$. CPU time 2064h.
\vspace*{3em}

\begin{figure}[h]
\label{Fe11-linreg}
\vspace*{-3em}
\includegraphics[width=270pt]{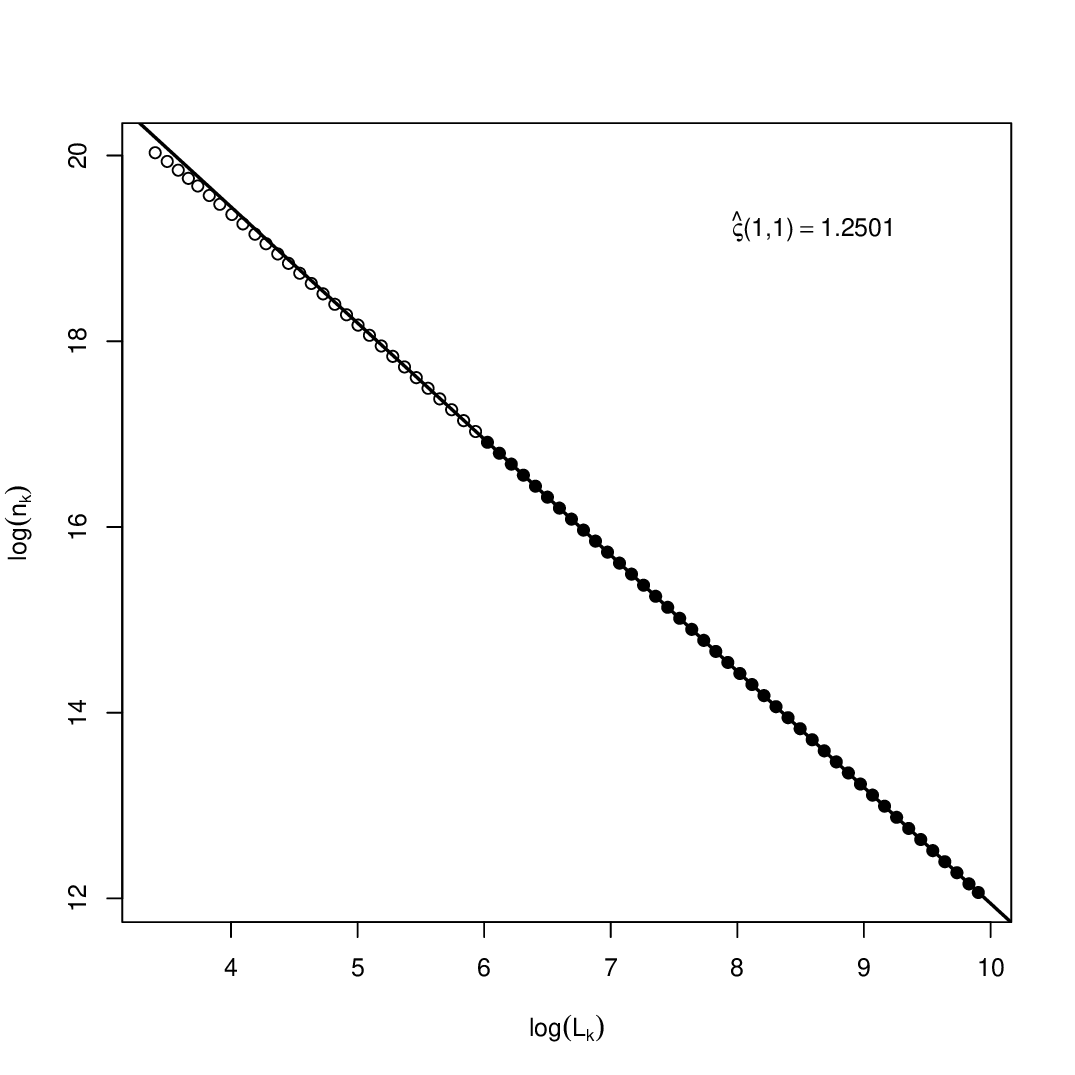}\vspace*{-2em}
\caption{Linear regression for the simulation of $\varsigma(1,1)$.}
\end{figure}

\begin{figure}[h]
\includegraphics[width=300pt]{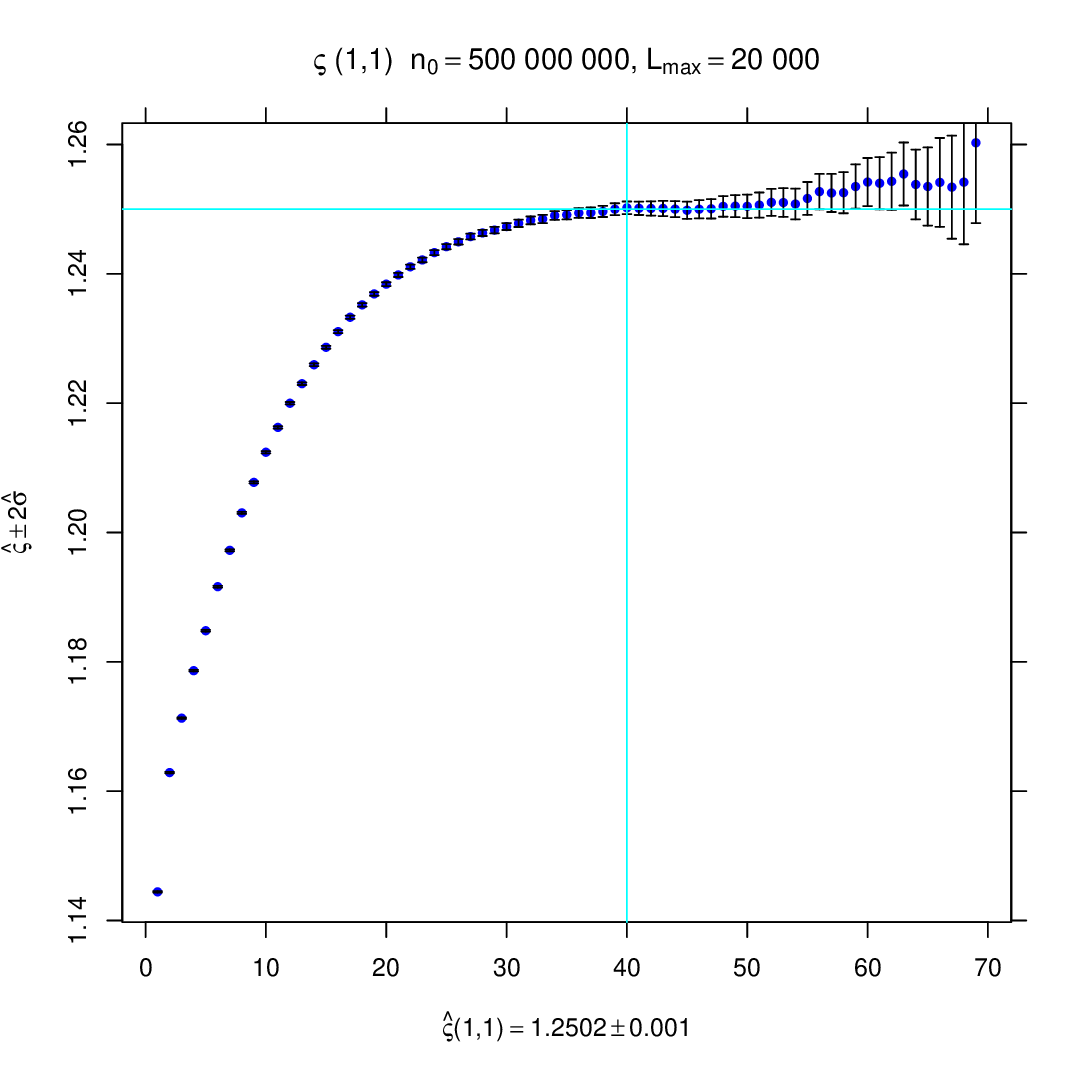}\vspace*{-0.1em}
\caption{Simulation for $\varsigma(1,1)$. The co-ordinate shows $\kmin$, the ordinate shows the corresponding $\hat\varsigma$ with error bars. The vertical line indicates $\kmin=40$ which we chose for our estimate of $\hat\varsigma$. The horizontal line shows the true value. }
\end{figure}
\clearpage

\subsubsection{Exponent $\varsigma(2,2)$.}
The exact value $\varsigma(2,2)=35/12=2.91666\ldots$ is known.
Also this simulation serves as a benchmark for our simulations.
\ftl
  30&  40000000000\\
   33&  27956276949\\
   36&  19934507109\\
   39&  14769670878\\
   42&  11270896745\\
   46&  8156016609\\
   50&  6108280379\\
   55&  4423365460\\
   60&  3315611015\\
   66&  2432628801\\
   72&  1842369408\\
   79&  1375726309\\
   86&  1056545535\\
   94&  803537797\\
  103&  607993657\\
 \ntl
  113&  459313243\\
  124&  347384944\\
  136&  263528000\\
  149&  200799711\\
  163&  153819037\\
  179&  116600065\\
  196&  89210485\\
  215&  67932955\\
  236&  51656232\\
  259&  39313221\\
  284&  30007400\\
  312&  22780638\\
  343&  17265563\\
  377&  13094893\\
  414&  9961095\\
 \ntl
  455&  7559087\\
  500&  5737717\\
  550&  4343548\\
  605&  3288311\\
  665&  2496057\\
  731&  1893876\\
  804&  1434709\\
  884&  1087314\\
  972&  823685\\
 1069&  624023\\
 1175&  473832\\
 1292&  359121\\
 1421&  271557\\
 1563&  205432\\
 1719&  155585\\
\ntl
 1890&  117893\\
 2079&  89442\\
 2286&  67757\\
 2514&  51314\\
 2765&  38803\\
 3041&  29341\\
 3345&  22363\\
 3679&  16949\\
 4046&  12813\\
 4450&  9738\\
 4895&  7339\\
 5384&  5571\\
 5922&  4218\\
 6514&  3229\\
 7165&  2477\\
\ntl
 7881&  1872\\
 8669&  1450\\
 9535&  1108\\
10488&  853\\
11536&  650\\
12689&  479\\
13957&  348\\
15352&  266\\
16887&  193\\
18575&  151\\
20000&  123\\
&\\&\\&\\&\\
\etl

Values used for the fit: $L_k=605\ldots20000$. CPU time 1879h.
\vspace*{-1em}
\begin{figure}[h]
\includegraphics[width=300pt]{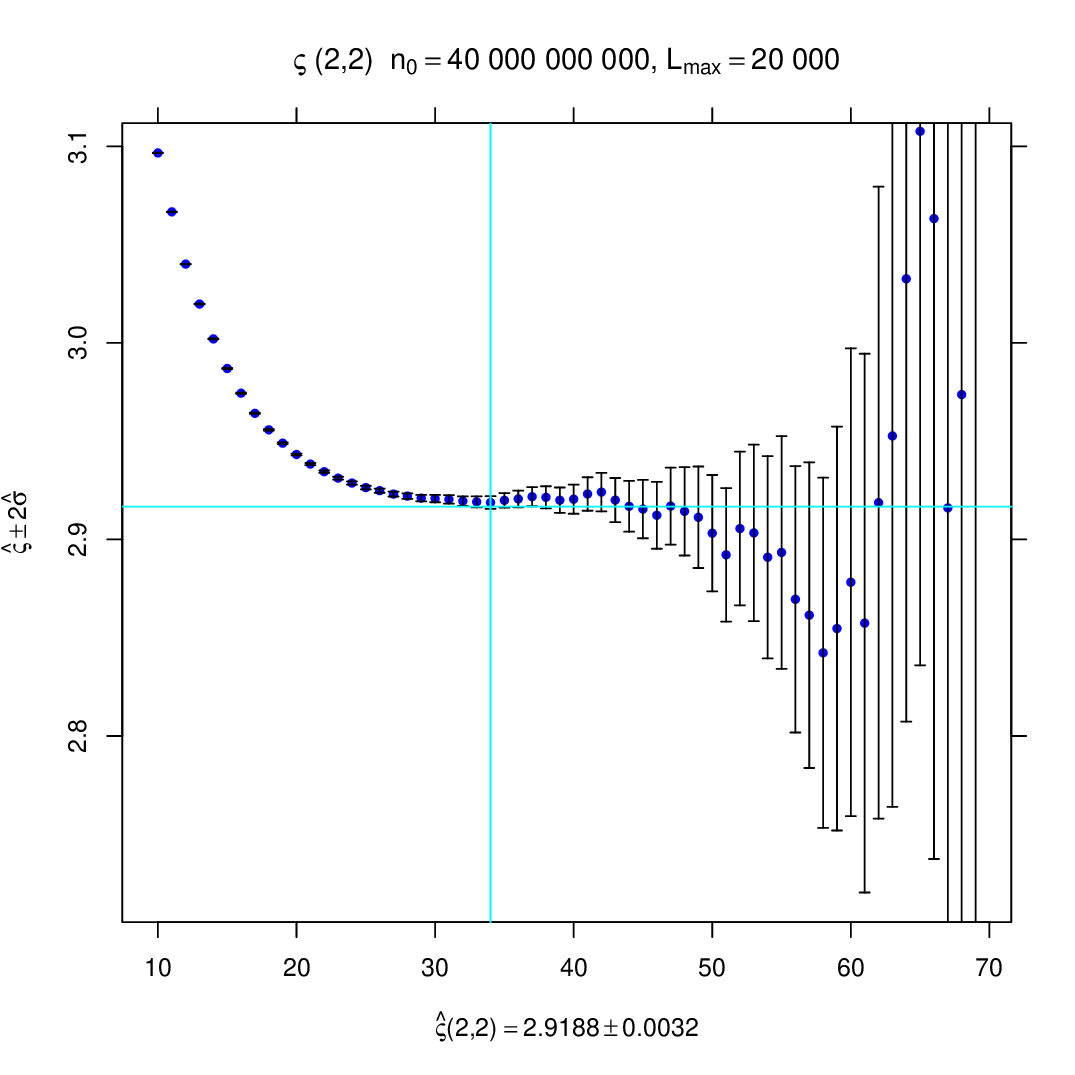}
\caption{Simulation for $\varsigma(2,2)$. The co-ordinate shows $\kmin$, the ordinate shows the corresponding $\hat\varsigma$ with error bars. The vertical line indicates $\kmin=34$ which we chose for our estimate of $\hat\varsigma$. The horizontal line shows the true value. }
\end{figure}
\clearpage

\subsubsection{Exponent $\varsigma(1,1,1)$.}
The exact value of $\varsigma(1,1,1)$ is unknown.

\ftl
   30&  60000000\\
   33&  59710616\\
   36&  58947709\\
   39&  57896946\\
   42&  56673833\\
   46&  54898618\\
   50&  53061195\\
   55&  50777396\\
   60&  48570208\\
   66&  46070175\\
   72&  43747356\\
   79&  41266693\\
   86&  39014779\\
   94&  36696389\\
  103&  34372986\\
  113&  32097711\\
  124&  29906035\\
  136&  27820941\\
\ntl
  149&  25859377\\
  163&  24028791\\
  179&  22226328\\
  196&  20584233\\
  215&  19009323\\
  236&  17526913\\
  259&  16147690\\
  284&  14873454\\
  312&  13666336\\
  343&  12537025\\
  377&  11494795\\
  414&  10540910\\
  455&  9652748\\
  500&  8835893\\
  550&  8076259\\
  605&  7376857\\
  665&  6740503\\
  731&  6155608\\
\ntl
  804&  5614962\\
  884&  5121770\\
  972&  4670981\\
 1069&  4256241\\
 1175&  3879722\\
 1292&  3534606\\
 1421&  3218775\\
 1563&  2930010\\
 1719&  2666485\\
 1890&  2426899\\
 2079&  2208165\\
 2286&  2009516\\
 2514&  1827541\\
 2765&  1661614\\
 3041&  1510468\\
 3345&  1372149\\
 3679&  1246493\\
 4046&  1132343\\
\ntl
 4450&  1028587\\
 4895&  934379\\
 5384&  848491\\
 5922&  770449\\
 6514&  699523\\
 7165&  635064\\
 7881&  576119\\
 8669&  523319\\
 9535&  474777\\
10488&  430885\\
11536&  391134\\
12689&  354964\\
13957&  321882\\
15352&  291736\\
16887&  264553\\
18575&  239803\\
20432&  217707\\
22475&  197364\\
\ntl
24722&  179046\\
27194&  162421\\
29913&  147273\\
32904&  133514\\
36194&  121333\\
39813&  109856\\
43794&  99544\\
48173&  90226\\
52990&  81910\\
58289&  74069\\
64117&  67148\\
70528&  60809\\
77580&  54981\\
80000&  53301\\
&\\&\\&\\&\\
\etl

Values used for the fit: $L_k=18575\ldots80\,000$. CPU time 8262h.
\begin{figure}[h]
\includegraphics[width=330pt]{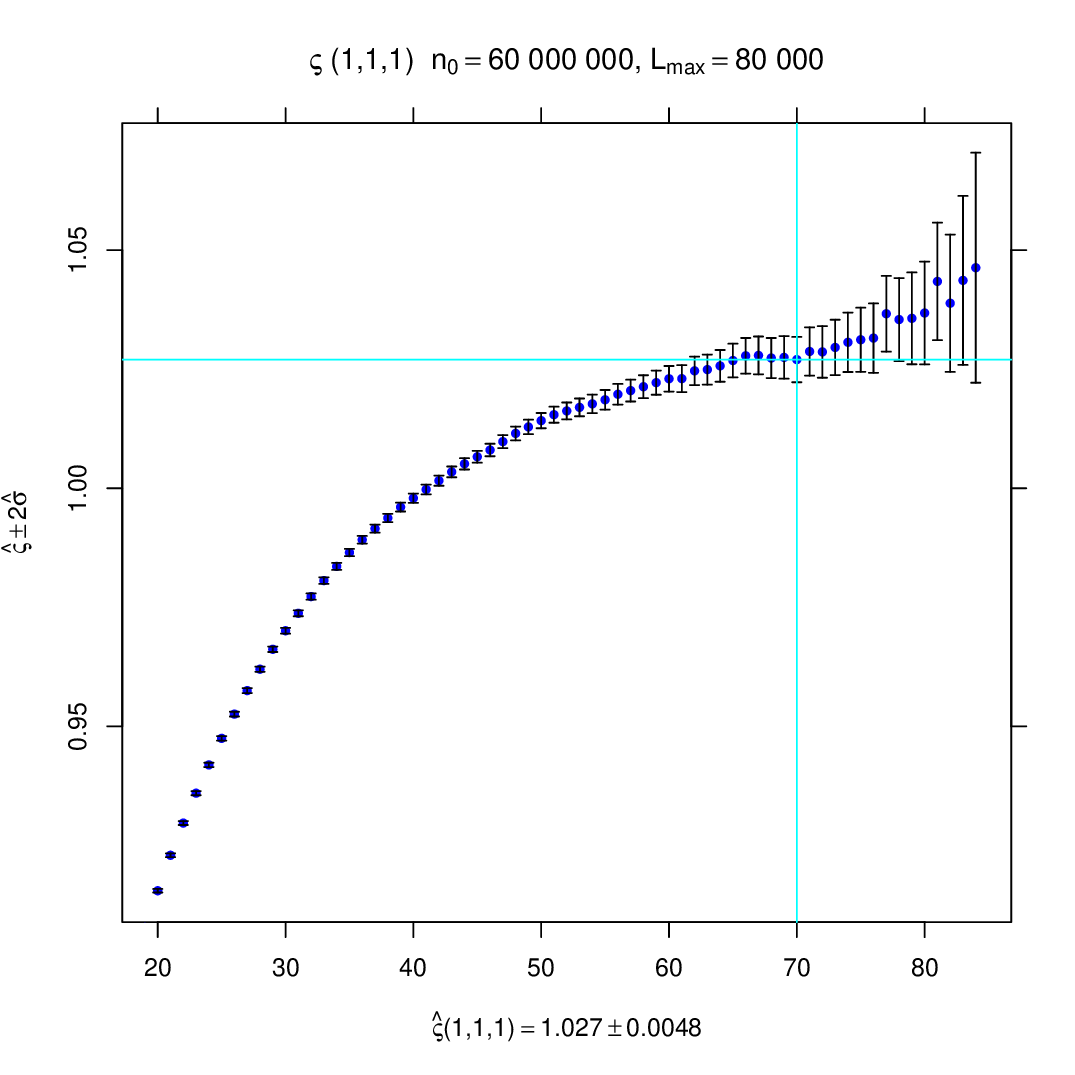}
\caption{Simulation for $\varsigma(1,1,1)$. The co-ordinate shows $\kmin$, the ordinate shows the corresponding $\hat\varsigma$ with error bars. The vertical line indicates $\kmin=70$ which we chose for our estimate of $\hat\varsigma$. The horizontal line shows the estimated value. }
\end{figure}
\subsubsection{Exponent $\varsigma(1,1,2)$.}
The exact value of $\varsigma(1,1,2)$ is unknown.

\ftl
  30&  200000000\\
   33&  198136199\\
   36&  193436538\\
   39&  187242650\\
   42&  180308394\\
   46&  170682110\\
   50&  161176218\\
   55&  149901066\\
   60&  139521065\\
   66&  128301060\\
   72&  118361321\\
   79&  108205442\\
   86&  99385254\\
   94&  90671582\\
  103&  82304628\\
  113&  74445095\\
  124&  67187208\\
  136&  60566796\\
\ntl
  149&  54585738\\
  163&  49210675\\
  179&  44119702\\
  196&  39647836\\
  215&  35523866\\
  236&  31776899\\
  259&  28415711\\
  284&  25417521\\
  312&  22668175\\
  343&  20191522\\
  377&  17981958\\
  414&  16028002\\
  455&  14267282\\
  500&  12694594\\
  550&  11279842\\
  605&  10020648\\
  665&  8909164\\
  731&  7917614\\
\ntl
  804&  7029965\\
  884&  6245336\\
  972&  5545792\\
 1069&  4923405\\
 1175&  4374033\\
 1292&  3885012\\
 1421&  3449618\\
 1563&  3062025\\
 1719&  2718548\\
 1890&  2415286\\
 2079&  2144282\\
 2286&  1904636\\
 2514&  1690316\\
 2765&  1499756\\
 3041&  1331441\\
 3345&  1181425\\
 3679&  1048523\\
 4046&  930691\\
\ntl
 4450&  826443\\
 4895&  733837\\
 5384&  651874\\
 5922&  578486\\
 6514&  513593\\
 7165&  456239\\
 7881&  405240\\
 8669&  359523\\
 9535&  319391\\
10488&  283792\\
11536&  251871\\
12689&  223767\\
13957&  198585\\
15352&  176146\\
16887&  156219\\
18575&  138744\\
20432&  123285\\
22475&  109424\\
\ntl
24722&  97235\\
27194&  86246\\
29913&  76472\\
32904&  67873\\
36194&  60371\\
39813&  53674\\
43794&  47628\\
48173&  42353\\
52990&  37627\\
58289&  33387\\
64117&  29608\\
70528&  26339\\
77580&  23407\\
80000&  22541\\
&\\&\\&\\&\\
\etl

Values used for the fit: $L_k=1069\ldots10\,000$. CPU time 5858h.
\begin{figure}[h]
\includegraphics[width=330pt]{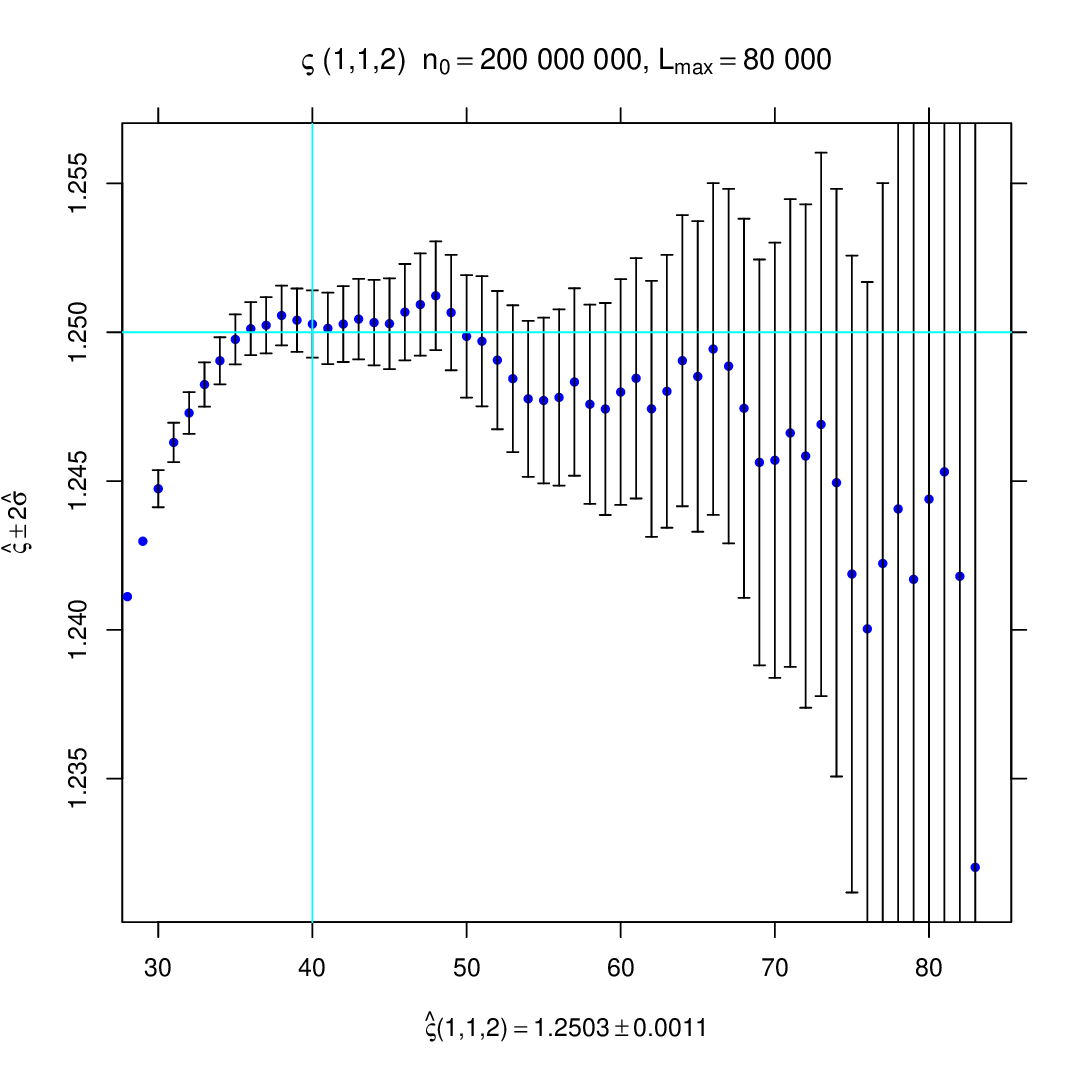}
\caption{Simulation for $\varsigma(1,1,2)$. The co-ordinate shows $\kmin$, the ordinate shows the corresponding $\hat\varsigma$ with error bars. The vertical line indicates $\kmin=40$ which we chose for our estimate of $\hat\varsigma$. The horizontal line shows the conjectured value $\varsigma(1,1,2)=\varsigma(1,1)=5/4$. }
\end{figure}
\clearpage
\subsubsection{Exponent $\varsigma(1,1,1,2)$.}
The exact value of $\varsigma(1,1,1,2)$ is unknown.
\ftl
30&  200000000\\
   33&  199921620\\
   36&  199482984\\
   39&  198565228\\
   42&  197176986\\
   46&  194683952\\
   50&  191624492\\
   55&  187237710\\
   60&  182471375\\
   66&  176514402\\
   72&  170515344\\
   79&  163642608\\
   86&  157023575\\
   94&  149860203\\
  103&  142343356\\
  113&  134666159\\
\ntl
  124&  126992806\\
  136&  119464557\\
  149&  112181604\\
  163&  105219661\\
  179&  98203442\\
  196&  91672706\\
  215&  85304312\\
  236&  79199999\\
  259&  73434697\\
  284&  68044645\\
  312&  62863201\\
  343&  57973159\\
  377&  53412940\\
  414&  49196512\\
  455&  45240841\\
  500&  41572496\\
\ntl
  550&  38132997\\
  605&  34954676\\
  665&  32040238\\
  731&  29348812\\
  804&  26850724\\
  884&  24560465\\
  972&  22451887\\
 1069&  20511329\\
 1175&  18739197\\
 1292&  17110256\\
 1421&  15612645\\
 1563&  14238301\\
 1719&  12983348\\
 1890&  11838167\\
 2079&  10786022\\
 2286&  9827571\\
\ntl
 2514&  8950110\\
 2765&  8149013\\
 3041&  7419539\\
 3345&  6752347\\
 3679&  6144676\\
 4046&  5589138\\
 4450&  5082774\\
 4895&  4621823\\
 5384&  4201858\\
 5922&  3819442\\
 6514&  3470994\\
 7165&  3155561\\
 7881&  2866913\\
 8669&  2604947\\
 9535&  2366048\\
10488&  2149715\\
\ntl
11536&  1953289\\
12689&  1773558\\
13957&  1609927\\
15352&  1461067\\
16887&  1326172\\
18575&  1203707\\
20432&  1092672\\
22475&  991657\\
24722&  900187\\
27194&  816464\\
29913&  740704\\
32904&  672302\\
36194&  609756\\
39813&  553485\\
40000&  550828\\
&\\\etl

Values used for the fit: $L_k=27\,194,  \ldots,40\,000  $. CPU time
35\,212h.
\begin{figure}[h]
\includegraphics[width=330pt]{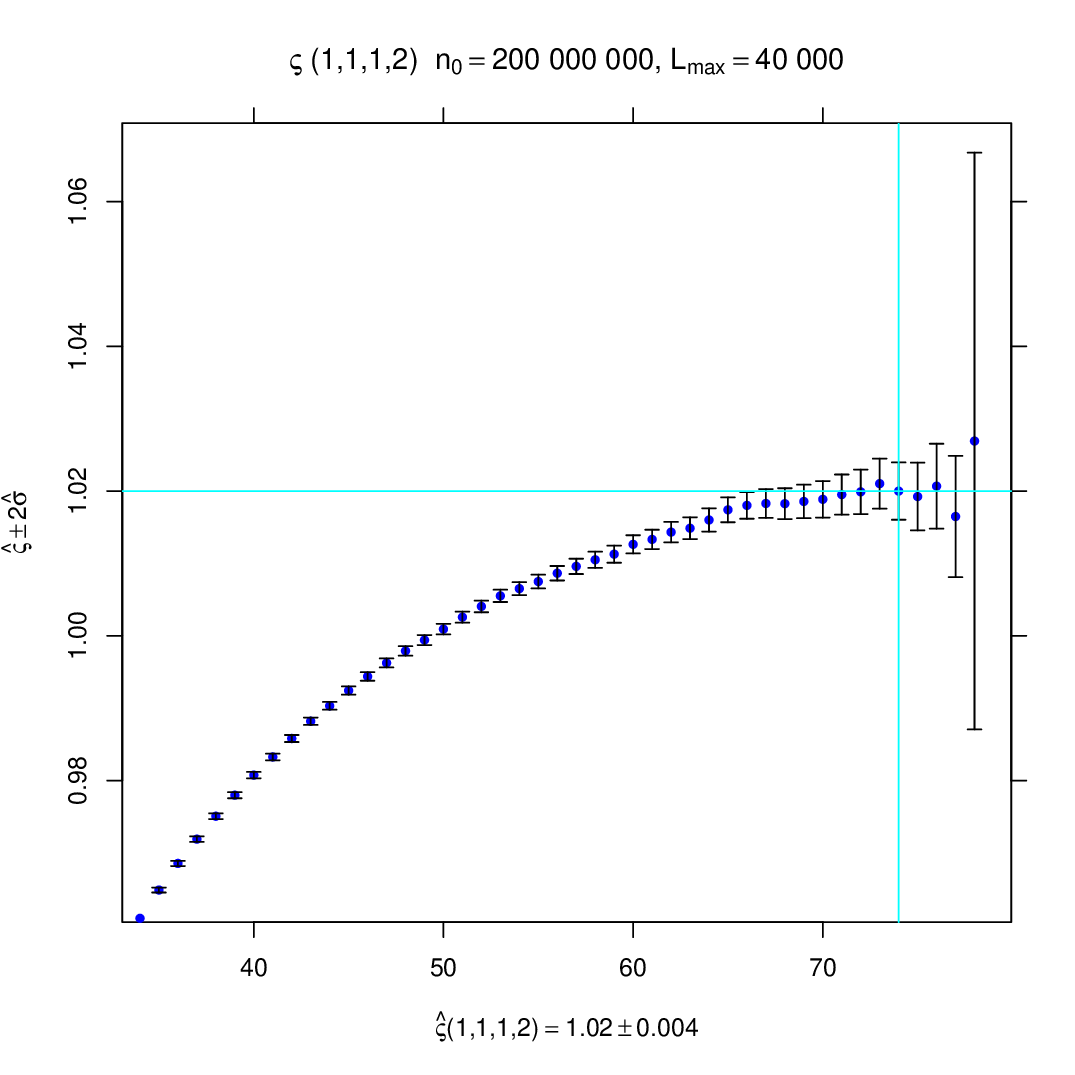}
\caption{Simulation for $\varsigma(1,1,1,2)$. The co-ordinate shows $\kmin$, the ordinate shows the corresponding $\hat\varsigma$ with error bars. The vertical line indicates $\kmin=74$ which we chose for our estimate of $\hat\varsigma$. The horizontal line shows the estimated value.}
\end{figure}

\clearpage
\subsubsection{Exponent $\varsigma(1,1,1,1)$.}
The exact value of $\varsigma(1,1,1,1)$ is unknown.
\ftl
  30&  20000000\\
   33&  19996035\\
   36&  19972855\\
   39&  19923361\\
   42&  19846358\\
   46&  19704559\\
   50&  19524417\\
   55&  19258012\\
   60&  18958581\\
   66&  18573948\\
   72&  18172653\\
   79&  17698995\\
   86&  17228707\\
   94&  16704112\\
  103&  16136193\\
  113&  15537061\\
  124&  14922694\\
  136&  14299200\\
\ntl
  149&  13677483\\
  163&  13066040\\
  179&  12433477\\
  196&  11828110\\
  215&  11220785\\
  236&  10622843\\
  259&  10040743\\
  284&  9483483\\
  312&  8935078\\
  343&  8402393\\
  377&  7892203\\
  414&  7410458\\
  455&  6946326\\
  500&  6504855\\
  550&  6081417\\
  605&  5681289\\
  665&  5304532\\
  731&  4949662\\
\ntl
  804&  4613179\\
  884&  4297577\\
  972&  4001670\\
 1069&  3722426\\
 1175&  3461745\\
 1292&  3217081\\
 1421&  2987877\\
 1563&  2773980\\
 1719&  2573328\\
 1890&  2386906\\
 2079&  2213006\\
 2286&  2051046\\
 2514&  1899964\\
 2765&  1758768\\
 3041&  1628269\\
 3345&  1506688\\
 3679&  1393480\\
 4046&  1289159\\
\ntl
 4450&  1191616\\
 4895&  1101052\\
 5384&  1017194\\
 5922&  939586\\
 6514&  867259\\
 7165&  801006\\
 7881&  739646\\
 8669&  682392\\
 9535&  629677\\
10488&  580582\\
11536&  535316\\
12689&  493266\\
13957&  454757\\
15352&  419131\\
16887&  386266\\
18575&  356146\\
20432&  328229\\
22475&  302420\\
\ntl
24722&  278878\\
27194&  256535\\
29913&  236157\\
32904&  217434\\
36194&  200369\\
39813&  184289\\
43794&  169567\\
48173&  156118\\
52990&  143655\\
58289&  132091\\
64117&  121391\\
70528&  111643\\
77580&  102601\\
80000&  99860\\
&\\
&\\
&\\
&\\
\etl

Values used for the fit: $L_k=39\,813\ldots80\,000$. CPU time 18\,262h.
\begin{figure}[h]
\includegraphics[width=330pt]{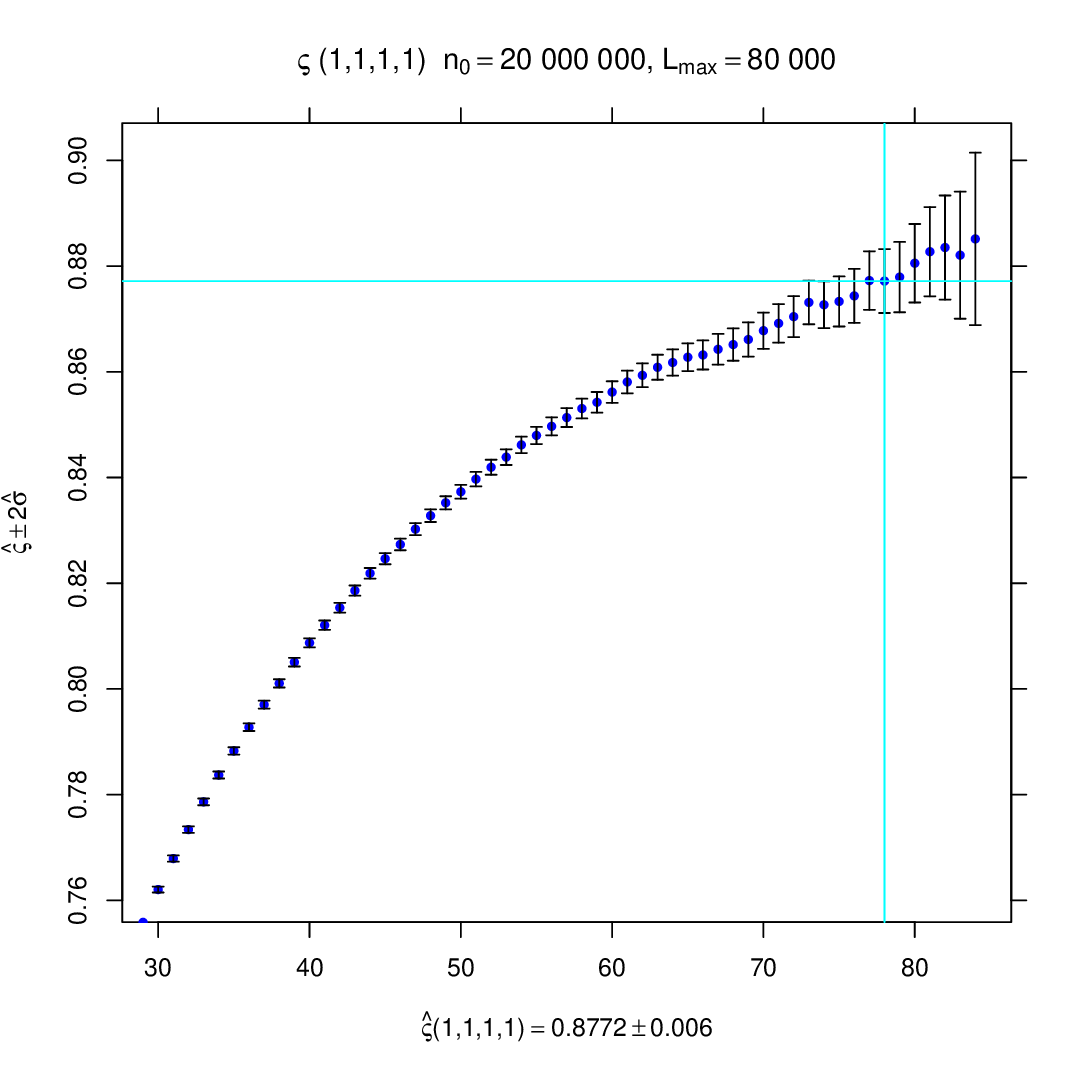}
\caption{Simulation for $\varsigma(1,1,1,1)$. The co-ordinate shows $\kmin$, the ordinate shows the corresponding $\hat\varsigma$ with error bars. The vertical line indicates $\kmin=78$ which we chose for our estimate of $\hat\varsigma$. The horizontal line shows the estimated value.}
\end{figure}

\clearpage
\subsubsection{Exponent $\varsigma(1,1,1,1,1)$.}
The exact value of $\varsigma(1,1,1,1,1)$ is unknown.
\ftl
  30&  744165\\
   33&  744158\\
   36&  744107\\
   39&  743886\\
   42&  743487\\
   46&  742538\\
   50&  741155\\
   55&  738765\\
   60&  735484\\
   66&  730711\\
   72&  725183\\
   79&  718008\\
   86&  710226\\
   94&  700782\\
  103&  689769\\
  113&  677384\\
\ntl
  124&  663493\\
  136&  648801\\
  149&  633473\\
  163&  617246\\
  179&  599664\\
  196&  581885\\
  215&  562961\\
  236&  543820\\
  259&  524462\\
  284&  505282\\
  312&  485535\\
  343&  465955\\
  377&  446300\\
  414&  426779\\
  455&  407418\\
  500&  388692\\
\ntl
  550&  370165\\
  605&  352289\\
  665&  335283\\
  731&  318821\\
  804&  302365\\
  884&  286416\\
  972&  271268\\
 1069&  256593\\
 1175&  242487\\
 1292&  229014\\
 1421&  216285\\
 1563&  204002\\
 1719&  192310\\
 1890&  181055\\
 2079&  170456\\
 2286&  160399\\
\ntl
 2514&  150828\\
 2765&  141592\\
 3041&  133025\\
 3345&  124845\\
 3679&  117257\\
 4046&  109961\\
 4450&  102996\\
 4895&  96520\\
 5384&  90434\\
 5922&  84699\\
 6514&  79371\\
 7165&  74357\\
 7881&  69601\\
 8669&  65234\\
 9535&  61002\\
10488&  56977\\
\ntl
11536&  53144\\
12689&  49724\\
13957&  46458\\
15352&  43355\\
16887&  40574\\
18575&  37983\\
20432&  35448\\
22475&  33083\\
24722&  30843\\
27194&  28666\\
29913&  26699\\
32904&  24955\\
36194&  23207\\
39813&  21616\\
40000&  21535\\
&\\
\etl

Values used for the fit: $L_k=27\,194\ldots40\,000$. CPU time 1147h.
\begin{figure}[h]
\includegraphics[width=330pt]{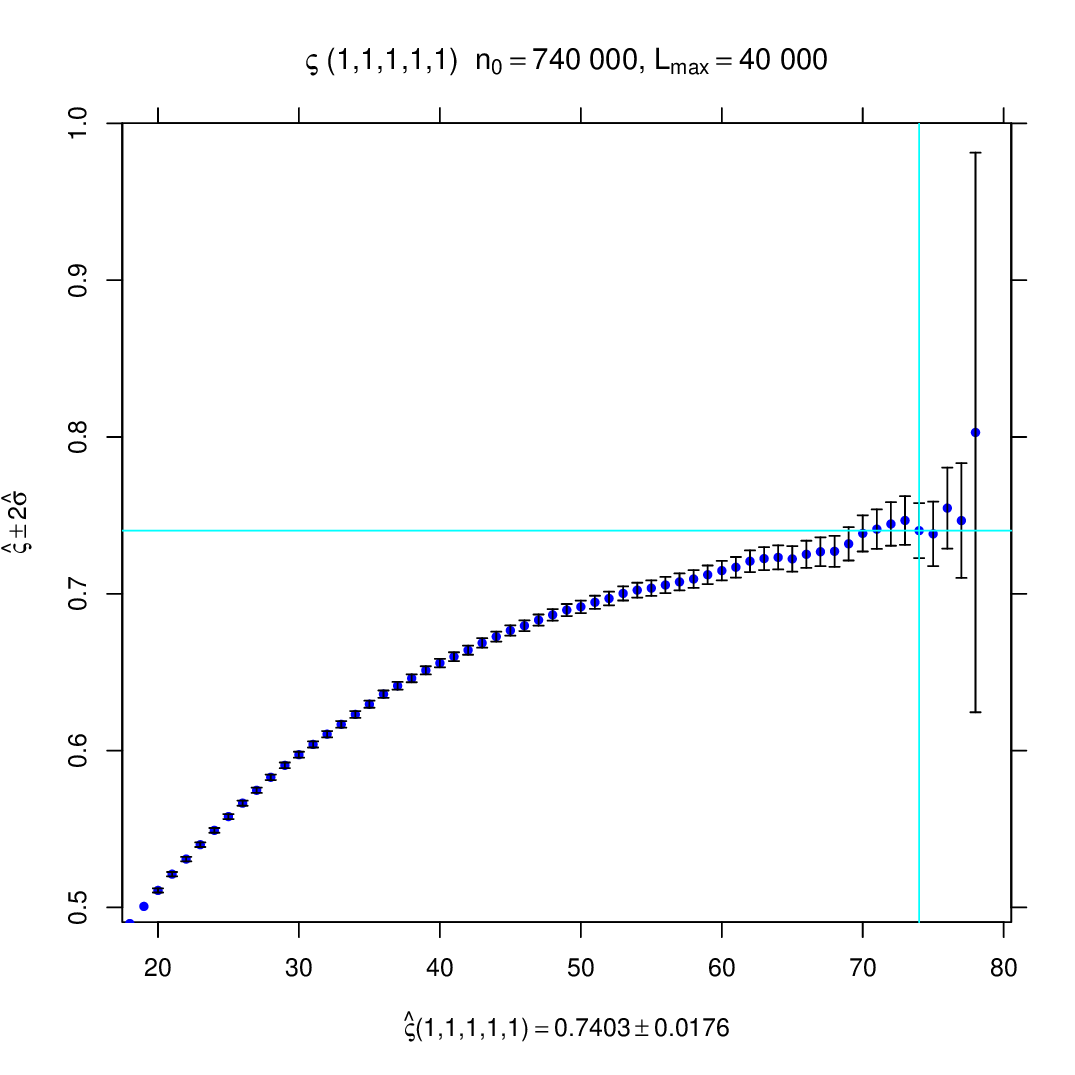}
\caption{Simulation for $\varsigma(1,1,1,1,1)$. The co-ordinate shows $\kmin$, the ordinate shows the corresponding $\hat\varsigma$ with error bars. The vertical line indicates $\kmin=74$ which we chose for our estimate of $\hat\varsigma$. The horizontal line shows the estimated value.}
\end{figure}

\subsubsection{Exponent $\varsigma(1,3,3)$.}
The exact value of $\varsigma(1,3,3)$ is unknown.
As it turns out that $\varsigma(1,3,3)>2$, we have performed simulations according to our scheme 2. That is, we have generated $n$ master samples of random walk paths that reach the boundary of the $L_1$-box (here $L_1=10000$). For each such master sample $i$ we have run $m=1000$ trials and have counted the fraction $\hat p_i$ of trials where the paths reached the boundary of the $L_2$-box (with $L_2=2L_1$). As $n=18\,100$ we cannot give the complete data set $p_1,\ldots,p_n$ but rather give the empirical mean and the standard deviation of $\hat p$
$$\hat p=0.155202983425414,\qquad \hat\sigma_p=0.000536918044881792.$$
From this we compute
$$\hat\varsigma(1,3,3)=2.6877718045551$$
with standard deviation
$$\hat\sigma=0.00499094143436367.$$
We conjecture that
$$\varsigma(1,3,3)=\varsigma(1,3)=\frac{13+\sqrt{73}}{8}=2.693000\ldots$$
We conclude with a histogram of the values $p_i$:
\begin{figure}[h]
\includegraphics[width=330pt]{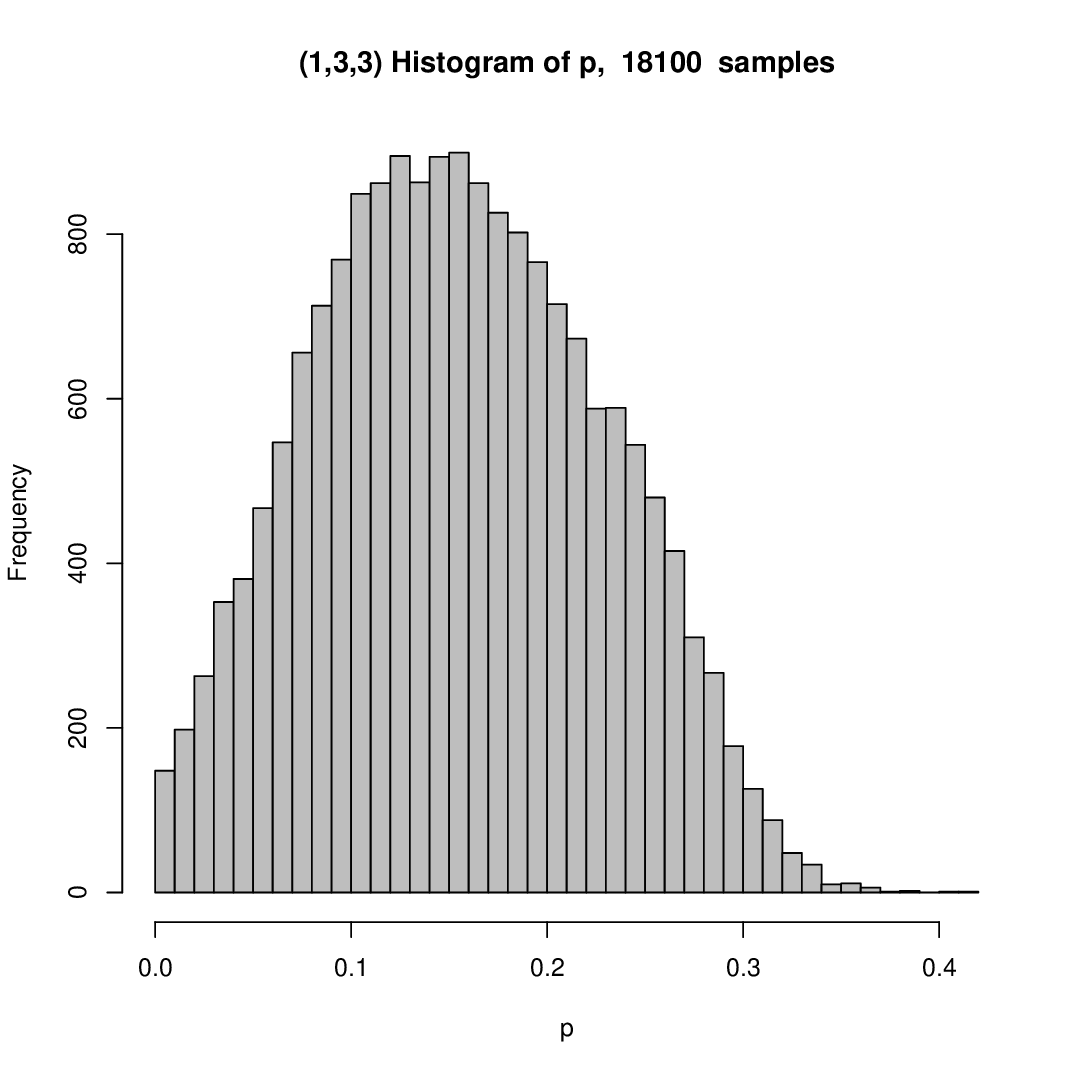}
\caption{Histogram of the values $p_i$ for  $\varsigma(1,3,3)$.}
\end{figure}
\subsubsection{Exponent $\varsigma(2,2,2)$.}
The exact value of $\varsigma(2,2,2)$ is unknown.
We have performed a simulation with the second scheme with $N=16000$, $n=1000$, $L_1=10000$, $L_2=20000$. Mean and standard deviation are
$$\hat p=0.1449495,\qquad \hat\sigma_p=0.000497221297799643.$$
From this we compute
$$\hat\varsigma(2,2,2)=2.78637773802317$$
with standard deviation
$$\hat\sigma=0.00494888703003405.$$
\subsubsection{Exponent $\varsigma(2,2,3)$.}
The exact value of $\varsigma(2,2,3)$ is unknown. We conjecture
$$\varsigma(2,2,3)=\varsigma(2,2)=\frac{35}{12}=2.916666\ldots$$

We have performed a simulation with the second scheme with $N=23000$, $n=1000$, $L_1=10000$, $L_2=20000$. Mean and standard deviation are
$$\hat p=0.130559,\qquad \hat\sigma_p=0.000444444142417374.$$
From this we compute
$$\hat\varsigma(2,2,3)=2.93722618256156$$
with standard deviation
$$\hat\sigma=0.00491116935805033.$$
\subsubsection{Exponent $\varsigma(2,3,3)$.}
The exact value of $\varsigma(2,3,3)$ is unknown.
We conjecture
$$\varsigma(2,3,3)=\varsigma(2,3)=\frac{47+5\sqrt{73}}{24}=3.738334113\ldots$$
We have performed a simulation with the second scheme with $N=1000$, $n=1000$, $L_1=10000$, $L_2=20000$. Mean and standard deviation are
$$\hat p=0.073458,\qquad \hat\sigma_p=0.00144828442088002.$$
From this we compute
$$\hat\varsigma(2,3,3)=3.76693657262376$$
with standard deviation
$$\hat\sigma=0.0284439101500224.$$
This simulation was particularly time consuming (179\,543h CPU time) as the actual value of $\varsigma(2,3,3)$ is rather large and it thus takes a tremendous amount of time to generate each master sample.
\subsubsection{Exponent $\varsigma(2,2,2,2)$.}
The exact value of $\varsigma(2,2,2,2)$ is unknown.
We have performed a simulation with the second scheme with $N=16000$, $n=1000$, $L_1=10000$, $L_2=20000$. Mean and standard deviation are
$$\hat p=0.157732125,\qquad \hat\sigma_p=0.000521232849038418.$$
From this we compute
$$\hat\varsigma(2,2,2,2)=2.66445157389522$$
with standard deviation
$$\hat\sigma=0.00476745017196814.$$



{\footnotesize
\begin{thebibliography}{BLP89}
\bibitem[BLP89]{BLP89}
{\sc K. Burdzy, G. Lawler} and {\sc T. Polaski.}
\newblock On the critical exponent for random walk intersections.
\newblock {\em J. Statist. Phys.} 56, 1-12 (1989).
\smallskip
\bibitem[CN06]{CN06}
{\sc F. Camia} and {\sc C. M. Newman.}
\newblock Two-dimensional critical percolation: the full scaling limit.
\newblock {\em Comm. Math. Phys.}, 268, 1--38 (2006).
\smallskip
\bibitem[CN07]{CN07}
{\sc F. Camia} and {\sc C. M. Newman.}
\newblock Critical percolation exploration path and SLE$_6$: a proof of convergence.
\newblock {\em Probab. Theory Related Fields} 139, 473--519 (2007).
\smallskip
\bibitem[Ca92]{Ca92}
{\sc J.~Cardy.}
\newblock Critical percolation in finite geometries.
\newblock {\em J.~Phys. A}, 25, L201-L206 (1992).
\smallskip
\bibitem[Du98]{Du98}
{\sc B.~Duplantier.}
\newblock Random walks and quantum gravity in two dimensions.
\newblock {\em Phys. Rev. Lett.}, 81, 5489-5492 (1998).
\smallskip
\bibitem[DK88]{DK88}
{\sc B.~Duplantier} and {\sc K.-H.~Kwon.}
\newblock Conformal invariance and intersections of random walks.
\newblock {\em Phys. Rev. Lett.}, 61, 2514-1517 (1988).
\smallskip
\bibitem[DS08]{DS08}
{\sc B.~Duplantier} and {\sc S.~Sheffield.}
\newblock Liouville Quantum Gravity and KPZ.
\newblock arXiv:0808.1560 (August 2008)
\smallskip
\bibitem[KM05]{KM05}
{\sc A. Klenke} and {\sc P.~M\"orters.}
\newblock The multifractal spectrum of Brownian intersection local time.
\newblock {\em Ann. Probab.} 33, 1255-1301 (2005).
\smallskip
\bibitem[Kn05]{Kn05}
{\sc D.~E.~Knuth.}
\newblock The art of computer programming.
\newblock Vol. {2}, 3rd edition. 
Addison-Wesley, Boston MA (2005).%
\smallskip
\bibitem[La89]{La89}
{\sc G.~F.~Lawler.}
\newblock Intersections of random walks with random sets.
\newblock {\em Israel J. Math.} 65, 113-132 (1989).
\smallskip
\bibitem[La91]{La91}
{\sc G.~F.~Lawler.}
\newblock Intersections of random walks.
\newblock Birkh\"auser, Boston MA (1991).
\smallskip
\bibitem[La95]{La95}
{\sc G.~F.~Lawler.}
\newblock Nonintersecting planar Brownian motions.
\newblock {\em Math. Phys. El. Journal}, 1, Paper~4, pp 1-35 (1995).%
\smallskip
\bibitem[La96]{La96}
{\sc G.~F.~Lawler.}
\newblock Hausdorff dimension of cut points for Brownian motion.
\newblock {\em El. Journal Probab.}, 1, Paper~2, pp 1-20 (1996).
\smallskip
\bibitem[LS01a]{LS01a}
{\sc G.~F.~Lawler, O.~Schramm} and {\sc W.~Werner.}
\newblock Values of Brownian intersection exponents I: Half-plane exponents.
\newblock {\em Acta Math.}, 187, 237-273 (2001).
\smallskip
\bibitem[LS01b]{LS01b}
{\sc G.~F.~Lawler, O.~Schramm} and {\sc W.~Werner.}
\newblock Values of Brownian intersection exponents II: Plane exponents.
\newblock {\em Acta Math.}, 187 , 275--308 (2001).
\smallskip
\bibitem[LS02]{LS02}
{\sc G.~F.~Lawler, O.~Schramm} and {\sc W.~Werner.}
\newblock Values of Brownian intersection exponents III: Two-sided exponents.
\newblock {\em Ann. Inst. Henri Poincar\'e}, 38, 109-123 (2002).
\smallskip
\bibitem[LR78]{LR78}
{\sc P. L. Leath} and {\sc G R Reich.}
\newblock Scaling form for percolation cluster sizes and perimeters.
\newblock {\em J. Phys. C} 11, 4017--4036 (1978).
\smallskip
\bibitem[Le83]{Le83}
{\sc E.~L.~Lehmann.}
\newblock Theory of point estimation.
\newblock {\em Wiley, New York} (1983).
\smallskip
\bibitem[LS90]{LS90}
{\sc B. Li} and {\sc A. Sokal.}
\newblock High-precision Monte Carlo test of the conformal-invariance predictions for two-dimensional mutually avoiding walks.
\newblock {\em J. Statist. Phys.} 61, 723--748 (1990).%
\smallskip
\bibitem[MP09]{MP09}
{\sc P.~M\"orters} and {\sc Y.Peres.}
\newblock Brownian motion.
\newblock {\em Cambridge University Press} (2009).
\smallskip
\bibitem[MS09]{MS09}
{\sc P.~M\"orters} and {\sc N.-R. Shieh.}
\newblock The exact packing measure of Brownian double points.
\newblock {\em Probab. Theory Related Fields} 143, 113--136 (2009).
\smallskip
\bibitem[RV08]{RV08}
{\sc R.~Rhodes} and {\sc V.~Vargas.}
\newblock KPZ formula for log-infinitely divisible multifractal random measures.
\newblock arXiv:0807.1036 (July 2008).
\smallskip
\bibitem[SD87]{SD87}
{\sc B.~Saleur} and {\sc B.~Duplantier.}
\newblock Exact determination of the percolation hull exponent in two dimensions.
\newblock {\em Phys. Rev. Lett.}, 58, 2325 - 2328 (1987).
\smallskip
\bibitem[Sm01]{Sm01}
{\sc S. Smirnov.}
\newblock Critical percolation in the plane: conformal invariance, Cardy's formula, scaling limits.
\newblock {\em C. R. Acad. Sci. Paris Sr. I Math.} 333, 239--244 (2001).
\smallskip
\bibitem[Vo84]{Vo84}
{\sc R.F. Voss.}
\newblock The fractal dimension of percolation cluster hulls.
\newblock {\em J. Phys. A} 17, L373-L377 (1984).
\smallskip
\end{thebibliography}}
\end{document}